\documentclass[11pt]{article}
\usepackage{amsfonts}
\usepackage{color}
\usepackage{graphics}
\usepackage{indentfirst}
\usepackage{cite}
\usepackage{latexsym}
\usepackage{amsmath}
\usepackage{amssymb}
\usepackage[dvips]{epsfig}
\usepackage{amscd}

\usepackage{leqno} 
\newtheorem{theorem}{Theorem}[section]
\newtheorem{remark}{Remark}[section]

\newtheorem{definition}{Definition}[section]
\newtheorem{lemma}[theorem]{Lemma}

\newtheorem{proposition}[theorem]{Proposition}

\newcommand{\no}{\nonumber\\}
\newcommand{\n}{\varrho}
\newcommand{\nn}{\varrho_0^R}

\newcommand{\ti}{\tilde}
\newcommand{\mr}{\mathbb{R}}

\def\pf{{\it Proof.}  }
\renewcommand{\div}{ {\rm div }  }

\newcommand{\uu}{\mathbf{u}}
\newcommand{\hh}{\mathbf{H}}
\newcommand{\na}{\nabla }
\newcommand{\vp}{\varphi }
\newcommand{\xx}{\mathbf{x}}
\newcommand{\pa}{\partial}

\newcommand{\bl}{\begin{lemma}}
\newcommand{\el}{\end{lemma}}
\newcommand{\ga}{\gamma}

\newcommand{\ve}{\varepsilon}
\newcommand{\la}{\label}

\newcommand{\om}{\Omega}

\newcommand{\bn}{\begin{eqnarray}}
\newcommand{\en}{\end{eqnarray}}
\newcommand{\bnn}{\begin{eqnarray*}}
\newcommand{\enn}{\end{eqnarray*}}

\newcommand{\ben}{\begin{enumerate}}
\newcommand{\een}{\end{enumerate}}

\newcommand{\be}{\begin{equation}}
\newcommand{\ee}{\end{equation}}

\def\p{\partial}
\def\norm[#1]#2{\|#2\|_{#1}}

\def\o{\omega}

\def\O{B_R}
\allowdisplaybreaks
\makeatletter      
\@addtoreset{equation}{section}
\makeatother
\begin{document}
\title{On classical solutions to the Cauchy problem of the 2D compressible non-resistive MHD equations with vacuum}
\date{}

\author{ Mingtao Chen$^{\lowercase{a}}$,   Aibin Zang$^{\lowercase{b}\thanks{Corresponding author}}$}
\maketitle
\begin{center}
$^a$ School of Mathematics and Statistics, Shandong University, \\
Weihai, Weihai 264209, P.R. China.\\
Email: mtchen@sdu.edu.cn

 $^b$  Center of Applied Mathematics, Yichun University, Yichun, Jiangxi, 336000, P.R.China\\
Email: 204124@jxycu.edu.cn 
\end{center}
%

%

\begin{abstract}
In this paper, we investigate the Cauchy problem of the compressible non-resistive MHD on $\mr^2$ with vacuum as far field density. We prove that the 2D Cauchy problem has a unique local strong solution provided the initial density and magnetic field decay not too slow at infinity. Furthermore, if the initial data satisfies some additional regularity and compatibility conditions, the strong solution becomes a classical one. Additionally, we establish a blowup criterion for the 2D compressible non-resistive MHD depending solely on the density and magnetic fields.
\end{abstract}

\textbf{Keywords}: 2D compressible non-resistive MHD equations; vacuum; classical solutions; blowup criterion

\section{Introduction}
In this paper, we study the two-dimensional (2D) compressible non-resistive magnetohydrodynamics (MHD) equations which read as follows:
\be\la{1.1}
\n_t+\div(\n \uu)=0,
\ee
\be\la{1.2}
(\n \uu)_t+\div(\n \uu\otimes \uu)+\nabla P(\n)=\mu \triangle \uu+(\mu+\lambda)\nabla\div \uu+(\na\times \hh)\times\hh,
\ee
\be\la{1.3}
\hh_t+\uu\cdot\na\hh+\hh\div\uu=\hh\cdot\na\uu,\quad \div\hh=0,
\ee
with the initial condition
\be\la{1.4}
(\n, \uu, \hh)(0, \xx)=(\n_0, \uu_0, \hh_0)(\xx),\quad \xx\in\mr^2,
\ee
and far field behavior (in some weak sense)
\be\la{1.5}
\uu(t, \xx)\rightarrow \mathbf{0},\quad \n(t, \xx)\rightarrow 0,\quad \hh(t, \xx)\rightarrow 0\quad {\rm as}\ |\xx|\rightarrow\infty, \quad {\rm for}\ t\geq0.
\ee
Here $\n=\n(t,\xx)$, $\uu(t, \xx)=(\uu_1, \uu_2)(t, \xx)$ and $\hh(t, \xx)=(\hh_1, \hh_2)(t, \xx)$ represent the unknown density, velocity and magnetic field of the fluid, respectively. The pressure $P(\n)$ is given by
\be\la{1.6}
P(\n)=A\n^\ga,
\ee
where $\ga>1$ is the adiabatic exponent, $A>0$ is a constant. The viscosity coefficients $\mu$ and $\lambda$ satisfy the following physical restrictions
\be\la{1.7}
\mu>0,\quad \mu+\lambda\geq0.
\ee

Magnetohydrodynamics studies the motion of electrically conducting media in the presence of a magnetic field. The dynamic motion of fluid and the magnetic field interact strongly with each other, so the hydrodynamic and electrodynamic effects are complicated either from the physical viewpoint, or from the mathematical consideration, see \cite{ca, ku} and references therein. In this paper, we restrict ourselves to \eqref{1.1}--\eqref{1.3} called the viscous and non-resistive magnetohydrodynamic equations, which means that the conducting fluids with a very high conductivity, such as, ideal conductors (c.f. \cite{cha, fre}). Particularly, the magnetic equation \eqref{1.3} implies that in a highly conducting fluid the magnetic field lines move along exactly with the fluid, rather than simply diffusing out. This type of behavior is physically expressed as that the magnetic field lines are frozen into the fluid. For more details of physical background, we refer the readers to \cite{ca, ku, cha, fre, xz} and references therein.

Now, we briefly recall some results concerning with the multi-dimensional compressible/incompressible MHD related with the present paper. First, if there is no electromagnetic effect, that is $\hh=0$, the MHD system reduces to the classical compressible/incompressible Navier-Stokes equations, which have been discussed by many mathematicians, please see \cite{cho, lions96, lions98, xin, fei, hlx12, hoff95, huang2, tem} and references therein. Next, if there exists electromagnetic effect with resistivity term, which also have been discussed by many researchers, we refer the reader to \cite{huang1, lxzsi, lv1, lv2, wang} and references therein.

However, as far as we known, the non-resistive MHD equations have not been thoroughly studied. As $\n\equiv$ constant, i.e., incompressible non-resistive MHD, Jiu et al. \cite{jiu} proved the local existence of solutions in 2D for initial data in $H^s$ for integer $s\geq3$, and soon after that Ren et al. \cite{ren} and Lin et al. \cite{lin} have established the existence of global-in-time solutions for initial data sufficiently close to certain equilibrium solutions in two spatial dimensions. Almost the same time, Fefferman et al. \cite{feff1} established local existence result for initial data in $(\uu_0, \hh_0)\in H^s$ with $s>d/2$ for $d=2, 3$; in order to assume less regularity on $\uu_0$ than that of $\hh_0$, due to the existence of the diffusive term in the momentum equations, Chemin et al. \cite{che} proved the local existence in Besov spaces when $\uu_0\in B_{2, 1}^{d/2-1}(\mr^d)$ and $\hh_0\in B_{2, 1}^{d/2}(\mr^d)$; and more recently, Fefferman et al. \cite{feff2} presented an inspiring local-in-time existence and uniqueness solutions in nearly optimal Sobolev space in $\mr^d$ ($d=2, 3$) for $\hh_0\in H^s(\mr^d)$ and $\uu_0\in H^{s-1+\ve}(\mr^d)$ with $s>d/2$ and $0<\ve<1$. Here we also want to mention Fan et al. \cite{zhou1} established the global-in-time existence of smooth solutions of 2D generalized MHD system with fractional diffusion $(-\triangle)^\alpha\uu$ $(0<\alpha<1/2)$.

Let's come back to the compressible non-resistive MHD equations. Fan et al. \cite{fan} and Li et al. \cite{lidd} independently proved the existence of unique local strong solution in 3D. After that, Xu et al. \cite{xz} established a blowup criterion to explain the mechanism of blow-up and the structure of possible singularities of strong solutions for the compressible non-resistive MHD equations, due to the lack of the global existence of strong solution to \eqref{1.1}--\eqref{1.3} with large initial data. Very recently, Zhu \cite{zhu} proved existence of unique local classical solution with regular initial data (almost the same compatibility conditions as that in \cite{fan}) and improved the blowup criterions obtained in \cite{xz}.

When the far field density is vacuum (particularly, the initial density may have compact support) in 2D, the methods  successfully applied to
3D in \cite{zhu, fan, lidd, xz} are not valid here, since the $L^p$-norm of $\uu$ could not bound in terms of $\|\sqrt\n\uu\|_{L^2(\mr^2)}$ and $\|\na\uu\|_{L^2(\mr^2)}$. Therefore, it is still unknown whether the well-posedness of strong/classical solutions to the compressible non-resistive MHD equations in 2D exist or not, even the local ones. In this paper, we want to answer parts of these questions.

In this section, for $1\leq r\leq \infty$, we denote the standard Lebesgue and Sobolev spaces as follows:
$$
L^r=L^r(\mr^2), \quad W^{s, r}=W^{s, r}(\mr^2), \quad H^s=W^{s, 2}.
$$

Now, we wish to define precisely what we mean by strong solutions.
\begin{definition}
If all derivatives involved in \eqref{1.1}--\eqref{1.3} for $(\n, \uu, \hh)$ are regular distributions, and equations \eqref{1.1}--\eqref{1.3} hold almost everywhere in $(0, T) \times \mr^2$, then $(\n, \uu, \hh)$ is called a strong solution to \eqref{1.1}--\eqref{1.3}.
\end{definition}

\begin{theorem}\la{thm1}
Let $\eta_0$ be a positive constant and
\be\la{1.8}
\bar\xx\triangleq (e+|\xx|^2)^{1/2}\ln^{1+\eta_0}(e+|\xx|^2).
\ee
For constants $q>2$ and $a>1$, assume that the initial data $(\n_0, \uu_0, \hh_0)$ satisfy that
\be\la{1.9}
\n_0\geq0, \ \n_0\bar\xx^a\in L^1\cap H^1\cap W^{1, q},\ \hh_0 \bar\xx^a \in  H^1\cap W^{1, q},\ \na\uu_0\in L^2, \sqrt{\n_0}\uu_0\in L^2.
\ee
Then there exists a positive time $T_0>0$ such that the problem \eqref{1.1}--\eqref{1.7} has a unique strong solution $(\n, \uu, \hh)$ on $(0, T_0]\times \mr^2$ satisfying that
\be\la{1.10}
\left\{
\begin{array}{lll}
\n\in C([0, T_0]; L^1\cap H^1\cap W^{1, q}), \ \n\bar\xx^a\in L^\infty (0, T_0; L^1\cap H^1\cap W^{1, q}),\\
\hh\in C([0, T_0]; H^1\cap W^{1, q}), \ \hh\bar\xx^a\in L^\infty(0, T_0; H^1\cap W^{1, q}),\\
\sqrt\n\uu,\ \na\uu, \ \bar\xx^{-1}\uu, \ \sqrt{t}\sqrt\n\uu_t\in L^\infty (0, T_0; L^2),\\
\na\uu\in L^2(0, T_0; H^1)\cap L^{(q+1)/q}(0, T_0; W^{1, q}), \ \sqrt{t} \na\uu\in L^2(0, T_0; W^{1, q}), \\
\sqrt\n \uu_t, \sqrt t\na\uu_t, \sqrt t\bar\xx^{-1}\uu_t\in L^2(\mr^2\times (0, T_0)),
\end{array}
\right.
\ee
and that
\be\la{1.11}
\inf_{0\leq t\leq T_0}\int_{B_N}\n(t, \xx)d\xx\geq \frac14\int_{\mr^2} \n_0(\xx)d\xx,
\ee
for some constant $N>0$ and $B_N\triangleq\{\xx\in\mr^2||\xx|<N\}$.
\end{theorem}

Moreover, if the initial data $(\n_0, \uu_0, \hh_0)$ satisfies some additional regularity and compatibility conditions, the local strong solution $(\n, \uu, \hh)$ obtained in Theorem \ref{thm1} becomes a classical one, that is,
\begin{theorem}\la{thm2}
In addition to \eqref{1.9}, suppose that
\be\la{1.12}
\left\{
\begin{array}{lll}
\na^2\n_0,\ \na^2 P(\n_0), \ \na^2\hh_0\in L^2\cap L^q,\\
\bar\xx^{\delta_0}\na^2\n_0, \ \bar\xx^{\delta_0}\na^2 P(\n_0), \ \bar\xx^{\delta_0}\na^2\hh_0, \ \na^2\uu_0\in L^2,
\end{array}
\right.
\ee
for some constant $\delta_0\in (0, 1)$. Moreover, assume that the following compatibility conditions hold for some $\mathbf{g}\in L^2$,
\be\la{1.13}
-\mu\triangle\uu_0-(\mu+\lambda)\na\div\uu_0+\na P(\n_0)-(\na\times \hh_0)\times\hh_0=\n_0^{1/2}\mathbf{g}.
\ee
Then, in addition to \eqref{1.10} and \eqref{1.11}, the strong solution $(\n, \uu, \hh)$ obtained in Theorem \ref{thm1} satisfies
\be\la{1.14}
\left\{
\begin{array}{lll}
\na^2\n,\ \na^2P(\n), \ \na^2\hh\in C([0, T_0]; L^2\cap L^q),\\
\bar\xx^{\delta_0}\na^2\n, \ \bar\xx^{\delta_0}\na^2 P(\n), \ \bar\xx^{\delta_0}\na^2\hh \in L^\infty(0, T_0; L^2),\\
\na^2\uu, \ \sqrt\n\uu_t, \ \sqrt t\na\uu_t,\ \sqrt t\bar\xx^{-1}\uu_t, \ t\sqrt\n \uu_{tt}, \ t\na^2\uu_t\in L^\infty(0, T_0; L^2),\\
t\na^3\uu\in L^\infty(0, T_0; L^2\cap L^q), \ \na\uu_t, \ \bar\xx^{-1} \uu_t, \ t\na\uu_{tt},\ t\bar\xx^{-1}\uu_{tt}\in L^2(0, T_0; L^2),\\
t\na^2(\n\uu)\in L^\infty(0, T_0; L^{(q+2)/2}).
\end{array}
\right.
\ee
\end{theorem}

There exists no global classical solution even without the effect of magnetic fields, due to \cite{xin} (c.f. \cite{yxin}). So, one naturally wonders: in finite time, what kinds of singularities will form, or what is the main mechanism of possible breakdown of smooth solutions for the 2D compressible non-resistive MHD equations? There are two main results \cite{xz, zhu} concerning blowup criteria for strong/classical solutions to the 3D compressible non-resistive MHD equations, which is similar as \cite{huang2} obtained for strong solutions to compressible Navier-Stokes equations. As it is well-known, partial differential equations (PDEs) entailing many independent variables are harder than PDEs entailing few independent variables. Therefore, there's an interesting question to ask whether the blowup criterions \cite{xz, zhu} could be improved for the 2D compressible non-resistive MHD or not. Based on subtle estimates, our next main result in this paper answered this question positively for classical solutions, which can be shown as follows.

\begin{theorem}\la{thm3}
Assume that the initial data $(\n_0, \uu_0, \hh_0)$ satisfies \eqref{1.9}, \eqref{1.12} and the compatibility conditions \eqref{1.13}. Let $(\n, \uu, \hh)$ be classical solution to the Cauchy problem \eqref{1.1}--\eqref{1.7}. If $0<T^*<\infty$ is the maximal time of existence, then
\be\la{1.15}
\underset{T\rightarrow T^*}{\lim\sup}\left(\|\n\|_{L^\infty(0, T; L^\infty)}+\|\hh\|_{L^\infty(0, T; L^\infty)}+\|\na\hh\|_{L^2(0, T; L^2(\mr^2))}\right)=\infty.
\ee
\end{theorem}

A few remarks are in order:
\begin{remark}
To obtain the local existence and uniqueness of strong/classical solutions, in Theorem \ref{thm1} and \ref{thm2}, the compatibility conditions we need is much weaker than the ones used in \cite{fan}, similar as that of \cite{zhu} for 3-D compressible non-resistive MHD equations.
\end{remark}

\begin{remark}
When $\hh=0$, i.e., there is no magnetic field effect, \eqref{1.1}--\eqref{1.2} reduces to the compressible Navier-Stokes equations, and Theorems \ref{thm1} and \ref{thm2} are similar to
the results of Li et al. \cite{lliang}. Roughly speaking, we generalize the results of \cite{lliang} to the 2D compressible non-resistive MHD equations. Furthermore, Theorems \ref{thm1} and \ref{thm2} extend the corresponding three-dimensional results in \cite{fan} and \cite{zhu} to 2D problem.
\end{remark}

\begin{remark}
Recently, Wang \cite{wang} established a blowup criterion for 2D compressible MHD equations which solely depends on the uniform (in time) upper bound of the density $\n$, i.e.,
\be\la{a1.15}
\lim_{T\rightarrow T^*}\|\n\|_{L^\infty(0, T; L^\infty)}=\infty.
\ee
It is clear that the blowup criterion in \eqref{1.15} for compressible non-resistive MHD equations is more stronger than the one in \eqref{a1.15}. This is mainly due to the lack of resistivity term for $\hh$, which can improve the stability of the system. However, Theorem \ref{thm3} is indeed improved previous ones obtained in \cite{xz, zhu}.
\end{remark}

\begin{remark}
 Compared with the previous blowup criterion established in \cite{zhou2} by Zhou et al. for 2D incompressible MHD system with zero magnetic diffusivity, which depends only on the magnetic fields, precisely, 
$$
\lim_{T\rightarrow T^*}\|\na\hh\|_{L^1(0, T; BMO(\mr^2))},
$$
our result is relatively stronger due to the compressibility and the existence of vacuum. We also want to refer the readers to \cite{zhou3}, where they built a series of the blowup criterions for 2D generalized incompressible MHD system.
\end{remark}

\begin{remark}
Compared with the incompressible non-resistive MHD \cite{feff1, feff2, che}, we consider the compressible one. Moreover, we indeed assume less regularity on $\uu_0$ than that of $\hh_0$ due to the existence of the diffusive term in \eqref{1.2} (see \eqref{1.9} and \eqref{1.12}), although it may not be optimal.
\end{remark}

We now make some comments on the analysis of this paper. The key difficulty of studying such MHD equations lies in the non-resistivity of the magnetic equations. However, for the 2D case, when the far field density is vacuum, another main difficulty is to bound the $L^p(\mr^2)$-norm of $\uu$ compared with 3D case, that means the methods which have successfully used in \cite{fan, lidd, xz, zhu} can't be directly applied to our cases. Fortunately, previous results \cite{hoff95, lions96, lions98, lliang, lv1, lv2} have provided hope for solving this problem. Precisely, we use the weighted $L^p$-bounds for elements of Hilbert space, see \eqref{2.5} and \eqref{2.6} below. Furthermore, compared with \cite{hoff95, lions98, lliang}, for the 2D compressible non-resistive MHD equations, the strong coupling between the velocity vector field and the magnetic field, such as $\uu\cdot\na\hh$ and $(\na\times\hh)\times\hh$ (which does not appear in compressible Navier-Stokes equations), will bring us some new difficulties. We'll borrow some ideas from \cite{lv2}, that is, in order to control the term $\uu\cdot\na\hh$ and related terms, we need a spatial weighted mean estimate of $\hh$ and $\na\hh$ (see \eqref{3.33} and \eqref{4.1} below). Compared with \cite{lv2}, we have to face other difficulties caused by the lack of resistive term in magnetic equations. We manage to solve the problem, because \eqref{1.1} and \eqref{1.3} have the analogous structure from the mathematical view point. Therefore, the magnetic field could be treated in a similar manner as that used for density, although \eqref{1.3} is more complicated than that of \eqref{1.1}.

The rest of the paper is organized as follows: In Section 2, we collect some elementary facts and inequalities which will be needed in later analysis. Sections 3 and 4 are devoted to the a priori
estimates which are needed to obtain the local existence and uniqueness of strong/classical solutions. The main results Theorem \ref{thm1}, \ref{thm2} and \ref{thm3} are proved in Section 5.

\section{Preliminaries}

First, the following local existence theory on bounded balls, where the initial density is strictly away from vacuum, can be shown by similar arguments as in \cite{fan, lidd, zhu}.
\bl\la{l2.1}
For $R>0$ and $B_R=\{\xx\in\mr^2||\xx|<R\}$, assume that $(\n_0, \uu_0, \hh_0)$ satisfies
\be\la{2.1}
(\n_0, \uu_0, \hh_0)\in H^3(B_R),\ \inf_{\xx\in B_R}\n_0(\xx)>0.
\ee
Then there exist a small time $T_R>0$ and a unique classical solution $(\n, \uu, \hh)$ to the following initial-boundary-value problem
\be\la{2.2}
\left\{
\begin{array}{lll}
\n_t+\div(\n\uu)=0,\\
\n\uu_t+\n(\uu\cdot\na)\uu+\na P(\n)-\mu\triangle\uu-(\mu+\lambda)\na \div\uu-(\na\times\hh)\times\hh=0,\\
\hh_t+(\uu\cdot\na)\hh+\hh\div\uu=(\hh\cdot\na)\uu,\\
\uu =0,\xx\in\p B_R,\quad t>0,\\
(\n, \uu, \hh)(0, \xx)=(\n_0, \uu_0, \hh_0)(\xx), \quad\xx\in B_R,
\end{array}
\right.
\ee
on $B_R\times (0, T_R]$ such that
\be\la{2.3}
\left\{
\begin{array}{lll}
\n, \hh\in C([0, T_R]; H^3), \ \uu\in C([0, T_R]; H^3)\cap L^2(0, T_R; H^4), \\
\uu_t\in L^\infty(0, T_R; H^1)\cap L^2(0, T_R; H^2), \ \sqrt\n\uu_{tt}\in L^2(0, T_R; L^2), \\
\sqrt t\uu\in L^\infty(0, T_R; H^4),\ \sqrt t\uu_t\in L^\infty(0, T_R; H^2),\ \sqrt t\uu_{tt}\in L^2(0, T_R; H^1),\\
\sqrt t\sqrt\n \uu_{tt}\in L^\infty(0, T_R; L^2), \ t\uu_t\in L^\infty(0, B_R; H^3),\\
t\uu_{tt}\in L^\infty(0, T_R; H^1)\cap L^2(0, T_R; H^2), \ t\sqrt\n \uu_{ttt}\in L^2(0, T_R; L^2), \\
t^{3/2}\uu_{tt}\in L^\infty(0, T_R; H^2), \ t^{3/2}\uu_{ttt}\in L^2(0, T_R; H^1),\\
t^{3/2}\sqrt\n\uu_{ttt}\in L^\infty(0, T_R; L^2),
\end{array}
\right.
\ee
where we denote $L^2=L^2(B_R)$ and $H^k=H^k(B_R)$ for some positive integer $k$.
\el

Then, for either $\Omega=\mr^2$ or $\Omega=\O$ with $R\geq1$, the following weighted $L^p$-bounds for elements of Hilbert space $\ti D^{1, 2} (\Omega) \triangleq \{v\in H_{\rm loc}^1(\Omega)|\na v\in L^2(\Omega)\}$ will play a crucial role in our analysis, which can be found in \cite[Lemma 2.4]{lliang}.
\bl\la{l2.3}
Let $\bar{\mathbf{x}}$ be as in \eqref{1.9} and $\Omega=\mr^2$ or $\Omega=\O$ with $R\geq 1$. For $\ga>1$, suppose that $\n\in L^1(\Omega)\cap L^\ga(\Omega)$ is a non-negative function such that
\be\la{2.4}
M_1\leq\int_{B_{N_1}}\n d\xx,\quad \int_\om P(\n)d\xx\leq M_2,
\ee
for some positive constants $M_1$, $M_2$ and $N_1\geq1$ with $B_{N_1}\subset\Omega$. Then, there exists a positive constant $C$ depending only on $M_1, M_2, N_1, \ga$  and $\eta_0$ such that
\be\la{2.5}
\|\mathbf v\bar\xx^{-1}\|_{L^2}\leq C\|\sqrt\n\mathbf v\|_{L^2}+C\|\na \mathbf v\|_{L^2},
\ee
for any $\mathbf v\in\ti D^{1, 2}(\om)$.

Furthermore, for $\ve>0$ and $\eta>0$, there is a positive constant $C$ depending on $\ve, \eta, M_1, M_2, N_1, \ga$ and $\eta_0$ such that any $\mathbf v\in\ti D^{1, 2}(\om)$ satisfies
\be\la{2.6}
\|\mathbf{v}\bar\xx^{-\eta}\|_{L^{(2+\varepsilon)/\tilde{\eta}}} \leq C\|\sqrt\n\mathbf v\|_{L^2}+C\|\na \mathbf v\|_{L^2},
\ee
with $\tilde{\eta}=\min\{1,\eta\}$.
\el

Next, we consider the following Lam\'{e} system,
\be\la{2.8}
\left\{
\begin{array}{lll}
-\mu\triangle \uu-(\mu+\lambda)\na\div\uu=F, &{\rm in\ }B_R,\\
\uu=0, &{\rm on\ }\p B_R.
\end{array}
\right.
\ee

The proof of the following $L^p$-bound is similar to that of \cite[Lemma 12]{cho}.
\bl\la{l2.5}
Let $\uu\in W_0^{1, q}(B_R)$ be a weak solution of the system \eqref{2.8}, where $q> 1$. If $F\in W^{k, q}(B_R)$ for $k\geq 0$, then $\uu\in W^{k+2, q}(B_R)$ and
\be\la{2.9}
\|\uu\|_{W^{k+2, q}(B_R)}\leq C\|\mathbf F\|_{W^{k, q}(B_R)},
\ee
where $C$ independent of $R$.
\el

Then, for $\na^\perp\triangleq (-\p_2, \p_1)$, denoting the material derivative of $\dot f\triangleq f_t+\uu\cdot\na f$. We now state some elementary estimates which follow from Gagliardo-Nirenberg inequality and the standard $L^p$-estimate for the following elliptic system derived from the momentum equations in \eqref{1.2}:
$$
\triangle F=\div\left(\n\dot\uu-\hh\cdot\na\hh\right),\quad \mu\triangle \o=\na^\perp\cdot\left(\n\dot\uu- \hh\cdot\na\hh\right),
$$
where
$$
F\triangleq(2\mu+\lambda)\div\uu-P(\n)-\frac12|\hh|^2, \quad \o=\p_1\uu^2-\p_2\uu^1.
$$

The proofs of the following results are similar as that of \cite[Lemma 2.5]{lv1}.
\bl\la{l2.6} Let $(\n, \uu, \hh)$ be a classical solution of \eqref{1.1}--\eqref{1.7}. Then for $p\geq 2$ there exists a positive constant $C$ depending only on $p, \mu$ and $\lambda$ such that
\be\la{2.12}
\|\na F\|_{L^p(\mr^2)}+\|\na \o\|_{L^p(\mr^2)}\leq C\left(\|\n \dot\uu\|_{L^p(\mr^2)}+\||\hh||\na\hh|\|_{L^p(\mr^2)}\right),
\ee
\begin{align}\la{2.13}
&\|F\|_{L^p(\mr^2)}+\|\o\|_{L^p(\mr^2)}\\
\leq& C\left(\|\n \dot\uu\|_{L^2(\mr^2)}+\||\hh||\na\hh|\| _{L^2(\mr^2)}\right)^{(p-2)/p} \left(\|\na\uu\|_{L^2(\mr^2)}+\|P(\n)\| _{L^2(\mr^2)}+\|\hh\|_{L^4}^2\right)^{2/p},\nonumber
\end{align}
\begin{align}\la{2.14}
\|\na\uu\|_{L^p(\mr^2)}\leq &C\left(\|\n \dot\uu\|_{L^2(\mr^2)} +\||\hh||\na\hh|\|_{L^2(\mr^2)}\right)^{(p-2)/p}+C\|\hh\|_{L^{2p}}^2\\
&\cdot\left(\|\na\uu\|_{L^2(\mr^2)}+\|P(\n)\| _{L^2(\mr^2)}+\|\hh\|_{L^4}^2\right)^{2/p}+C\|P(\n)\|_{L^p(\mr^2)} .\nonumber
\end{align}
\el

Finally, the following Beale-Kato-Majda type inequality, which is similar as that of \cite[Lemma 2.3]{huang3} will be used later to estimate $\|\na\uu\|_{L^\infty}$ and $\|\na\n\|_{L^2\cap L^q}$ ($q>2$).
\bl\la{l2.6}
For $2<q<\infty$, there is a constant $C(q)$ such that the following estimate holds for all $\na\uu\in L^2(\mr^2)\cap D^{1, q}(\mr^2)$,
\begin{align}\la{2.15}
\|\na\uu\|_{L^\infty(\mr^2)}\leq& C\left(\|\div\uu\|_{L^\infty(\mr^2)} +\|\o\|_{L^\infty(\mr^2)}\right)\ln\left(e+\|\na^2\uu\|_{L^q(\mr^2)} \right)\\
&\quad+C\|\na\uu\|_{L^2(\mr^2)}+C.\nonumber
\end{align}\el

\section{ A  priori estimates (I)}
In this section and the next, for $p\in[1, \infty]$ and $k\geq0$, we denote
$$
\int f d\xx=\int_{B_R} fd\xx, \quad L^p=L^p(B_R), \quad W^{k, p}=W^{k, p}(B_R), \quad H^k=W^{k, 2}.
$$
Moreover, for $R>4N_0\geq 4$, we assume that the smooth triplet $(\n_0, \uu_0, \hh_0)$ satisfies, in addition to \eqref{2.1}, that
\be\la{3.1}
\frac12\leq \int_{B_{N_0}}\n_0(\xx)d\xx\leq \int_{B_R}\n_0(\xx)d\xx\leq \frac32.
\ee
It follows from Lemma \ref{2.1} that there exists some $T_R>0$ such that the initial-boundary-value problem \eqref{2.2} has a unique classical solution $(\n, \uu, \hh)$ on $[0, T_R]\times B_R$ satisfying \eqref{2.3}.
For $\bar\xx, \eta_0, a$ and $q$ as in Theorem \ref{thm1}, the main aim of this section is to derive the following key apriori estimate on $\phi(t)$, defined by
\be\la{3.2}
\phi(t)\triangleq 1+\|\sqrt\n\uu\|_{L^2}+\|\na\uu\|_{L^2}+\|\hh\|_{L^2} +\|\hh\bar\xx^a\|_{H^1\cap W^{1, q}}+\|\n\bar\xx^a\|_{L^1\cap H^1\cap W^{1, q}}.
\ee

\begin{proposition}\la{prop}
Assume that $(\n_0, \uu_0, \hh_0)$ satisfies \eqref{2.1} and \eqref{3.1}. Let $(\n, \uu, \hh)$ be the solution to the initial-boundary-value problem \eqref{2.2} on $(0, T_R]\times B_R$ obtained by Lemma \ref{2.1}. Then there exist positive constants $T_0$ and $M$ both depending on $\mu, \lambda, \ga, q, a, \eta_0, N_0$ and $C_0$ such that
\be\la{3.3}
\sup_{0\leq t\leq T_0}\phi(t)+\int_0^{T_0}\left(\|\na^2 \uu\|_{L^q} ^{(q+1)/q}+t\|\na^2 \uu\|_{L^q}^2+\|\na^2 \uu\|_{L^2}^2\right)dt\leq M,
\ee
where
$$
C_0=\|\sqrt{\n_0}\uu_0\|_{L^2}+\|\na\uu_0\|_{L^2}+\|\hh_0\|_{L^2} +\|\hh_0\bar\xx^a\|_{H^1\cap W^{1, q}}+\|\n_0\bar\xx^a\| _{L^1\cap H^1\cap W^{1, q}}.
$$
\end{proposition}

The proof of Proposition \ref{prop} will be postponed at the end of this section. First, we start with the following energy estimate for $(\n, \uu, \hh)$ and preliminary $L^2$-bounds for $\na\uu$.
\bl\la{l3.1}
Let $(\n, \uu, \hh)$ be a smooth solution to the initial-boundary-value problem \eqref{2.2}. Then there exist a positive constant $\alpha=\alpha(\ga, q)>1$ and a $T_1=T_1(C_0, N_0)>0$ such that for all $t\in (0, T_1]$,
\begin{align}\la{3.4}
&\sup_{0\leq s\leq t}\left(\|\na\uu\|_{L^2}^2+\|\hh\|_{L^2}^2+ \|\sqrt\n\uu\|_{L^2}^2 +\|P(\n)\|_{L^1}\right)\\
&\qquad\quad+\int_0^t\left(\|\na\uu\|_{L^2}^2+\|\sqrt\n \uu_t\|_{L^2}^2\right)ds
\leq C+C\int_0^t\phi^\alpha(s) ds.\nonumber
\end{align}
\el
\pf First, multiplying \eqref{2.2}$_2$ and \eqref{2.2}$_3$ by $\uu$ and $\hh$, respectively, and integrating the resultant equalities over $B_R$, and summing them together, then integration by parts show that
\be\la{3.5}
\sup_{0\leq s\leq t}\left(\|\sqrt\n\uu\|_{L^2}^2+\|\hh\|_{L^2}^2 +\|P(\n)\|_{L^1}\right) +\int_0^t\|\na\uu\|_{L^2}^2ds\leq C.
\ee

Next, for $N>1$ and $\vp_N\in C_0^\infty(B_R)$ such that
\be\la{3.6}
0\leq\vp_N\leq 1, \ \vp_N(\xx)=1,\ {\rm if}\ |\xx|\leq N/2,\ \ |\na^k\vp_N|\leq CN^{-k} (k=1, 2),
\ee
then it follows from \eqref{3.1} and \eqref{3.5} that
\begin{align}\la{3.7}
\frac{d}{dt}\int\n\vp_{2N_0}d\xx=&\int\n\uu\cdot\na\vp_{2N_0}d\xx\\
\geq&-CN_0^{-1}\left(\int\n d\xx\right)^{1/2}\left( \int\n|\uu|^2 \right)^{1/2}\geq -\ti C(C_0, N_0),\nonumber
\end{align}
where in the last inequality we have used
$$
\int\n d\xx=\int\n_0 d\xx,
$$
due to \eqref{2.2}$_1$ and \eqref{2.2}$_4$. Integrating \eqref{3.7} over $(0, T_1)$ shows
\be\la{3.8}
\inf_{0\leq t\leq T_1}\int_{B_{2N_0}}\n d\xx\geq \inf_{0\leq t\leq T_1} \int\n\vp_{2N_0}d\xx\geq \int\n_0\vp_{2N_0}d\xx-\ti CT_1\geq1/4,
\ee
where $T_1\triangleq\min\{1, (4\ti C)^{-1}\}$. From now on, we will always suppose that $t\leq T_1$. The combination of \eqref{2.6}, \eqref{3.5} and \eqref{3.8} shows that for $\ve>0$ and $\eta>0$, every $\mathbf v\in \ti D^{1, 2}(B_R)$ satisfies
\be\la{3.9}
\|\mathbf v\bar\xx^{-\eta}\|_{L^{(2+\ve)/\ti\eta}}^2\leq C(\ve, \eta)\|\sqrt\n\mathbf v\|_{L^2}^2+C(\ve, \eta)\|\na\mathbf v\|_{L^2}^2,
\ee
with $\ti\eta=\min\{1, \eta\}$. In particular, we have
\be\la{3.10}
\|\n^\eta\uu\|_{L^{(2+\ve)/\ti\eta}}+\|\uu\bar\xx^{-\eta}\| _{L^{(2+\ve)/\ti\eta}}\leq C(\ve, \eta)\phi^{1+\eta}(t).
\ee

Next, multiplying \eqref{2.2}$_2$ by $\uu_t$ and integration by parts yield
\begin{align}\la{3.11}
&\frac{d}{dt}\left(\mu\|\na\uu\|_{L^2}^2+(\mu+\lambda) \|\div\uu\|_{L^2}^2\right)+\|\sqrt\n\uu_t\|_{L^2}^2\\
\leq&C\int \n|\uu|^2|\na\uu|^2d\xx+2\int P(\n)\div\uu_t d\xx+2\int (\hh\cdot\na)\hh \cdot\uu_t d\xx+\int |\hh|^2\div\uu_t d\xx. \nonumber
\end{align}
Now we estimate each term on the right-hand side of \eqref{3.11}. First, the Gagliardo-Nirenberg inequality implies that for all $p\in(2, +\infty)$,
\be\la{3.12}
\|\na\uu\|_{L^p}\leq C\|\na\uu\|_{L^2}^{2/p}\|\na\uu\|_{H^1}^{1-2/p} \leq C\phi(t)+C\phi(t)\|\na^2\uu\|_{L^2}^{1-2/p},
\ee
which together with \eqref{3.10} yields that for $\eta>0$ and $\ti\eta=\min\{1, \eta\}$,
\begin{align}\la{3.13}
\int\n^\eta|\uu|^2|\na\uu|^2d\xx\leq & C\|\n^{\eta/2}\uu\|_{L^{8/\ti\eta} }^2\|\na\uu\|_{L^{8/(4-\ti\eta)}}^2\\
\leq& C(\eta)\phi^{4+2\eta}(t)\left(1+\|\na^2\uu\|_{L^2}^{\ti\eta/2} \right)\no
\leq&C\phi^{\alpha(\eta)}(t)+\ve\phi^{-1}(t)\|\na^2\uu\|_{L^2}^2.\nonumber
\end{align}

Next, since $P(\n)$ satisfies
\be\la{3.14}
P_t(\n)+\div(P(\n)\uu)+(\ga-1)P(\n)\div\uu=0,
\ee
we deduce from \eqref{3.10}, \eqref{3.13} and the Sobolev inequality that
\begin{align}\la{3.15}
&2\int P(\n)\div\uu_td\xx\\
=&2\frac{d}{dt}\int P(\n)\div\uu d\xx-2\int P(\n) \uu\cdot\na \div\uu d\xx+2(\ga-1)\int P(\n)(\div\uu)^2d\xx\no
\leq &2\frac{d}{dt}\int P(\n)\div\uu d\xx+\ve\phi^{-1}(t) \|\na^2\uu\|_{L^2}^2+C(\ve)\phi^{\alpha}(t).\nonumber
\end{align}

Then, using integration by parts together with \eqref{2.2}$_3$, one obtains
\begin{align}\la{3.16}
2\int (\hh\cdot\na)\hh \cdot\uu_t d\xx=&-2\frac{d}{dt}\int (\hh\cdot\na ) \uu\cdot\hh d\xx+2\int (\hh_t\cdot\na)\uu\cdot\hh d\xx\\
&\quad+2\int (\hh\cdot\na )\uu\cdot\hh_t d\xx\no
=&-2\frac{d}{dt}\int (\hh\cdot\na ) \uu\cdot\hh d\xx+2\int ((\hh\cdot\na)\uu\cdot\na)\uu\cdot\hh d\xx\no
&-2\int ((\uu\cdot\na)\hh\cdot\na)\uu\cdot\hh d\xx-2\int (\hh\div\uu\cdot\na)\uu\cdot\hh d\xx\no
&+2\int (\hh\cdot\na )\uu\cdot(\hh\cdot\na)\uu d\xx-2\int (\hh\cdot\na )\uu\cdot(\uu\cdot\na)\hh d\xx\no
&-2\int (\hh\cdot\na )\uu\cdot\hh\div\uu d\xx.\nonumber
\end{align}
First, it is easy to check that
\begin{align*}
&\left|\int ((\hh\cdot\na)\uu\cdot\na)\uu\cdot\hh d\xx\right|+\left|\int (\hh\cdot\na )\uu\cdot(\hh\cdot\na)\uu d\xx\right|\no
+&\left|\int (\hh\div\uu\cdot\na)\uu\cdot\hh d\xx\right|+\left|\int (\hh\cdot\na )\uu\cdot\hh\div\uu d\xx\right|
\leq C\|\hh\|_{L^\infty}^2\|\na\uu\|_{L^2}^2.
\end{align*}
Next, H\"{o}lder inequality and Young inequality, together with \eqref{3.12}, yield
\begin{align*}
&\left|\int ((\uu\cdot\na)\hh\cdot\na)\uu\cdot\hh d\xx\right| +\left|\int (\hh\cdot\na )\uu\cdot(\uu\cdot\na)\hh d\xx\right|\no
\leq &C\int|\hh||\na\uu||\uu||\na\hh|d\xx\\
\leq &C\|\hh\|_{L^\infty}\|\uu\bar\xx^{-1}\|_{L^{2a}}\|\hh\bar\xx^a\| _{H^1}\|\na\uu\|_{L^{2a/(a-1)}}\no
\leq &C\phi^\alpha(t)+\ve\phi^{-1}(t)\|\na^2\uu\|_{L^2}^2.
\end{align*}
Substituting the above two estimates into \eqref{3.16} gives
\be\la{3.17}
2\int (\hh\cdot\na)\hh \cdot\uu_t d\xx\leq-2\frac{d}{dt}\int (\hh\cdot\na ) \uu\cdot\hh d\xx +C\phi^\alpha(t)+\ve\phi^{-1}(t)\|\na^2\uu\|_{L^2}^2.
\ee
Similarly,
\be\la{3.18}
\int |\hh|^2\div\uu_t d\xx\leq \frac{d}{dt}\int |\hh|^2\div\uu d\xx  +C\phi^\alpha(t)+\ve\phi^{-1}(t)\|\na^2\uu\|_{L^2}^2.
\ee

Inserting \eqref{3.13}, \eqref{3.15}, \eqref{3.17} and \eqref{3.18} into \eqref{3.11} shows
\begin{align}\la{3.19}
&\frac{d}{dt}\left(\mu\|\na\uu\|_{L^2}^2+(\mu+\lambda) \|\div\uu\|_{L^2}^2\right)\\
&-2\frac{d}{dt}\int \left(P(\n)\div\uu +\frac12|\hh|^2\div\uu d\xx-2 (\hh\cdot\na ) \uu\cdot\hh\right) d\xx+\|\sqrt\n\uu_t\|_{L^2}^2\no
\leq &C\phi^\alpha(t)+4\ve\phi^{-1}(t)\|\na^2\uu\|_{L^2}^2.\nonumber
\end{align}

To estimate the last term on the right-hand side of \eqref{3.19}, it follows from \eqref{2.9} that for $p\in [2, q]$,
\begin{align}\la{3.20}
\|\na^2\uu\|_{L^p}\leq & C\left(\|\n\uu_t\|_{L^p}+\|\n\uu\cdot\na\uu\|_{L^p}+\|\na P(\n)\|_{L^p}+\||\hh||\na \hh|\|_{L^p}\right),
\end{align}
which together with \eqref{3.12} and \eqref{3.13} yields
\begin{align}\la{3.21}
\|\na^2\uu\|_{L^2}\leq& C\phi^{1/2}(t)\|\sqrt\n\uu_t\|_{L^2} +C\|\n\uu\cdot\na\uu\|_{L^2} +C\phi^\alpha (t)\\
\leq &C\phi^{1/2}(t)\|\sqrt\n\uu_t\|_{L^2}+\frac12\|\na^2\uu\|_{L^2} +C\phi^\alpha (t).\nonumber
\end{align}
Substituting \eqref{3.21} into \eqref{3.19}, then integrating the resultant inequality over $(0, t)$, and choosing $\ve$ suitably small leads to
\begin{align}\la{3.22}
&\frac\mu2\|\na\uu\|_{L^2}^2+(\mu+\lambda)\|\div\uu\|_{L^2}^2+\int_0^t \|\sqrt\n\uu_t\|_{L^2}^2ds\\
\leq &C+C\|P(\n)\|_{L^2}^2+C\|\hh\|_{L^4}^4+C\int_0^t\phi^\alpha(s)ds\no
\leq &C+C\|\hh\|_{L^4}^4+C\int_0^t\phi^\alpha(s)ds,\nonumber
\end{align}
where in the last inequality we have used
\begin{align*}
\|P(\n)\|_{L^2}^2\leq & C\|P(\n_0)\|_{L^2}^2+C\int_0^t\|P(\n)\| _{L^\infty}^{3/2}\|P(\n)\|_{L^1}^{1/2}\|\na\uu\|_{L^2}ds\\
\leq& C+C\int_0^t\phi^\alpha(s)ds,
\end{align*}
due to \eqref{3.14}.

To estimate the second term on the right-hand side of \eqref{3.22}, multiplying \eqref{2.2}$_3$ by $4|\hh|^2\hh$ and integrating the resultant equality over $B_R$, we have
\begin{align*}
\frac{d}{dt}\|\hh\|_{L^4}^4\leq C\int|\na\uu||\hh|^4d\xx
\leq C\|\na\uu\|_{L^2}\|\hh\|_{L^2}\|\hh\|_{L^\infty}^3.\nonumber
\end{align*}
Integrating the above inequality over $(0, t)$, yields
\be\la{3.23}
\|\hh\|_{L^4}^4\leq C+\int_0^t\phi^\alpha(s)ds.
\ee
Putting \eqref{3.23} into \eqref{3.22}, together with \eqref{3.5}, leads to \eqref{3.4}. Therefore, we complete the proof of Lemma \ref{3.1}. \hfill $\Box$

\bl\la{l3.2}
Let $(\n, \uu, \hh)$ and $T_1$ be as in Lemma \ref{l3.1}. Then for all $t\in (0, T_1]$,
\be\la{3.24}
\sup_{0\leq s\leq t}s\|\sqrt\n\uu_t\|_{L^2}^2+\int_0^ts\|\na\uu_t\| _{L^2}^2ds\leq C\exp\left(C\int_0^t\phi^\alpha(s)ds\right).
\ee
\el
\pf Differentiating \eqref{2.2}$_2$ with respect to $t$ yields
\begin{align}\la{3.25}
&\n\uu_{tt}+\n\uu\cdot\na\uu_t-\mu\triangle\uu_t-(\mu+\lambda) \na\div\uu_t \\
=&-\n_t(\uu_t+\uu\cdot\na\uu)-\n\uu_t\cdot\na\uu-\na P_t(\n)-\frac12 \na|\hh|_t^2+\hh_t\cdot\na\hh+\hh\cdot\na\hh_t.\nonumber
\end{align}
Multiplying \eqref{3.25} by $\uu_t$ and integrating the resultant equation over $B_R$, we obtain
\begin{align}\la{3.26}
&\frac12\frac{d}{dt}\int\n|\uu_t|^2d\xx+\mu\int|\na\uu_t|^2d\xx+ (\mu+\lambda)\int|\div\uu_t|^2d\xx\\
=&-2\int\n\uu\cdot\na\uu_t\cdot\uu_t d\xx-\int\n\uu\cdot\na(\uu\cdot \na\uu\cdot\uu_t)d\xx\no
&-\int\n\uu_t\cdot\na\uu\cdot\uu_t d\xx+\int P_t(\n)\div\uu_td\xx+ \frac12\int |\hh|^2_t\div\uu_td\xx\no
&-\int(\hh\otimes\hh)_t:\na\uu_td\xx\no
\leq &C\int\n|\uu||\uu_t|\left(|\na\uu_t|+|\na\uu|^2+|\uu||\na^2\uu |\right)d\xx+C\int\n|\uu|^2|\na\uu||\na\uu_t|d\xx\no
&+C\int\n|\uu_t|^2|\na\uu|d\xx+C\int |P_t(\n)||\div\uu_t|d\xx+C \int |\hh||\hh_t||\na\uu_t|d\xx.\nonumber
\end{align}

We now estimate each term on the right-hand side of \eqref{3.26} as follows:

First, it follows from \eqref{3.2}, \eqref{3.5}, \eqref{3.9}, \eqref{3.10} and \eqref{3.12} that for $\ve\in (0, 1)$,
\begin{align}\la{3.27}
&\int\n|\uu||\uu_t|\left(|\na\uu_t|+|\na\uu|^2+|\uu||\na^2\uu |\right)d\xx\\
\leq &C\|\sqrt\n\uu\|_{L^6}\|\sqrt\n\uu_t\|_{L^2}^{1/2} \|\sqrt\n\uu_t\|_{L^6}^{1/2}\left(\|\na\uu_t\|_{L^2} +\|\na\uu\|_{L^4}^2\right)\no
&+C\|\n^{1/4}\uu\|_{L^{12}}^2\|\sqrt\n\uu_t\|_{L^2}^{1/2}  \|\sqrt\n\uu_t\|_{L^6}^{1/2}\|\na^2\uu\|_{L^2}\no
\leq &C\phi^\alpha(t)\|\sqrt\n\uu_t\|_{L^2}^{1/2}\left( \|\sqrt\n\uu_t\|_{L^2}^{1/2}+ \|\na\uu_t\|_{L^2}^{1/2}\right) \left(\|\na\uu_t\|_{L^2}+\|\na^2\uu\|_{L^2}+\phi(t)\right)\no
\leq&\ve\|\na\uu_t\|_{L^2}^2+C\phi^\alpha(t)\left( \|\na^2\uu\|_{L^2}^2+\|\sqrt\n\uu_t\|_{L^2}^2+1\right).\nonumber
\end{align}

Next, H\"{o}lder inequality together with \eqref{3.10} and \eqref{3.12} shows that
\begin{align}\la{3.28}
\int\n|\uu|^2|\na\uu||\na\uu_t|d\xx\leq & C\|\sqrt\n\uu\|_{L^8}^2\|\na \uu\|_{L^4}\|\na\uu_t\|_{L^2}\\
\leq &\ve\|\na\uu_t\|_{L^2}^2+C\left(\phi^\alpha (t)+\|\na^2\uu\|_{L^2}^2 \right).\nonumber
\end{align}

Then, H\"{o}lder inequality and \eqref{3.9} yield that
\begin{align}\la{3.29}
\int\n|\uu_t|^2|\na\uu|d\xx\leq &\|\na\uu\|_{L^2}\|\sqrt\n \uu_t\|_{L^6}^{3/2} \|\sqrt\n\uu_t\|_{L^2}^{1/2}\\
\leq &\ve\|\na\uu_t\|_{L^2}^2+C\phi^\alpha(t)\|\sqrt\n \uu_t\|_{L^2}^2.\nonumber
\end{align}

Next, it follows from \eqref{3.10} and \eqref{3.14} that
\begin{align}\la{3.30}
&\int|P_t(\n)||\div\uu_t|d\xx\\
\leq & C\int \left(P(\n)|\div\uu|+|\na P(\n)||\uu|\right) |\div\uu_t|d\xx\no
\leq & C\left(\|P(\n)\|_{L^\infty}\|\na\uu\|_{L^2}+\|\n\|_{L^\infty} ^{\ga-1} \|\n\bar\xx^a\|_{W^{1, q}}\|\uu\bar\xx^{-a}\|_{L^{2q/(q-2)}} \right)\|\na\uu_t\|_{L^2}\no
\leq &\ve\|\na\uu_t\|_{L^2}^2+C\phi^\alpha(t).\nonumber
\end{align}

Finally, H\"{o}lder inequality, \eqref{2.2}$_3$ and \eqref{3.10} give that
\begin{align}\la{3.31}
\int |\hh||\hh_t||\na\uu_t|d\xx
\leq &C\int \left(|\hh||\na\uu|+|\na\hh||\uu|\right)|\na\uu_t|d\xx\\
\leq &C\left(\|\hh\|_{L^\infty}\|\na\uu\|_{L^2}+\|\hh\bar\xx^a\|_{W^{1, q}} \|\uu\bar\xx^{-a}\|_{L^{2q/(q-2)}}\right)\|\na\uu_t\|_{L^2}\no
\leq &\ve\|\na\uu_t\|_{L^2}^2+C\phi^\alpha(t).\nonumber
\end{align}

Inserting \eqref{3.27}--\eqref{3.31} into \eqref{3.26} and choosing $\ve$ suitably small yield that
\begin{align}\la{3.32}
&\frac{d}{dt}\int\n|\uu_t|^2d\xx+\mu\int|\na\uu_t|^2d\xx+ (\mu+\lambda)\int|\div\uu_t|^2d\xx\\
\leq &C\phi^\alpha(t)\left(1+\|\sqrt\n\uu_t\|_{L^2}^2 +\|\na^2\uu\|_{L^2}^2 \right)\no
\leq &C\phi^\alpha(t)\|\sqrt\n\uu_t\|_{L^2}^2+C\phi^\alpha(t),\nonumber
\end{align}
where in the last inequality we have used \eqref{3.21}. Then, multiplying \eqref{3.32} by $t$, we finally obtain \eqref{3.24} after using Gronwall inequality and \eqref{3.4}. Therefore, we complete the proof of Lemma \ref{l3.2}. \hfill $\Box$

\bl\la{l3.3}
Let $(\n, \uu, \hh)$ and $T_1$ be as in Lemma \ref{l3.1}. Then for all $t\in (0, T_1]$,
\be\la{3.33}
\sup_{0\leq s\leq t}\left(\|\n\bar\xx^a\|_{L^1\cap H^1\cap W^{1, q}}+\|\hh\bar\xx^a\|_{H^1\cap W^{1, q}}\right) \leq \exp\left(C\exp\left(C\int_0^t\phi^\alpha(s)ds\right)\right).
\ee
\el
\pf First, multiplying \eqref{2.2}$_1$ by $\bar\xx^a$ and integrating the resultant equality over $B_R$, integration by parts and using \eqref{3.10}, we have
\begin{align*}
\frac{d}{dt}\int\n\bar\xx^ad\xx\leq &C\int \n |\uu|\bar\xx^{a-1} \ln^{1+\eta_0}(e+|\xx|^2)d\xx\\
\leq &C\|\n\bar\xx^{a-1+8/(8+a)}\|_{L^{(8+a)/(7+a)}}\|\uu \bar\xx^{-4/(8+a)}\|_{L^{8+a}}\\
\leq &C\phi^\alpha(t),
\end{align*}
which implies
\be\la{3.34}
\sup_{0\leq s\leq t}\|\n\bar\xx^a\|_{L^1}\leq C\exp\left(C \int_0^t\phi^\alpha(s)ds\right).
\ee

Next, it follows from the Sobolev inequality and \eqref{3.10} that for $0<\delta<1$,
\begin{align}\la{3.35}
\|\uu\bar\xx^{-\delta}\|_{L^\infty}\leq &C\left(\|\uu\bar\xx^{-\delta} \|_{L^{4/\delta}}+\|\na(\uu\bar\xx^{-\delta})\|_{L^3} \right)\\
\leq& C\left(\|\uu\bar\xx^{-\delta} \|_{L^{4/\delta}}+ \|\na\uu\|_{L^3}+\|\uu\bar\xx^{-\delta} \|_{L^{4/\delta}}\|\bar\xx^{-1}\na\bar\xx\|_{L^{12/(4-3\delta)}} \right)\no
\leq &C\left(\phi^\alpha (t) +\|\na^2\uu\|_{L^2}\right).\nonumber
\end{align}

Then, one derives from \eqref{2.2}$_3$ that $\ti\hh\triangleq \hh\bar\xx^a$ satisfies
\be\la{3.36}
\ti\hh_t+\uu\cdot\na \ti\hh-a\uu\cdot\na\ln\bar\xx\cdot\ti\hh +\ti\hh\div\uu =\ti\hh\cdot\na\uu,
\ee
which together with \eqref{3.35} and Gagliardo-Nirenberg inequality  shows that
\begin{align}\la{3.37}
\frac{d}{dt}\|\ti\hh\|_{L^2}\leq &C\left(\|\na\uu\|_{L^\infty} +\|\uu\cdot\na\ln\bar\xx\|_{L^\infty}\right)\|\ti\hh\|_{L^2}\\
\leq &C\left(\phi^\alpha(t)+\|\na^2\uu\|_{L^2\cap W^{1, q}}\right)\|\ti\hh\|_{L^2}.\nonumber
\end{align}
Moreover, \eqref{3.36} together with \eqref{3.35} and Gagliardo-Nirenberg inequality gives that for $p\in [2, q]$
\begin{align}\la{3.38}
&\frac{d}{dt}\|\na\ti\hh\|_{L^p}\\
\leq & C\left( 1+\|\na\uu\|_{L^\infty} +\|\uu\cdot\na\ln\bar\xx\|_{L^\infty}\right)\|\na\ti\hh\|_{L^p}\no
&+C\left(\||\na\uu||\na\ln\bar\xx|\|_{L^p}+\||\uu||\na^2\ln\bar\xx|\| _{L^p} +\|\na^2\uu\|_{L^p}\right)\|\ti\hh\|_{L^\infty}\no
\leq &C\left(\phi^\alpha(t)+\|\na\uu\|_{L^2\cap W^{1, q}}\right)\|\na\ti\hh\|_{L^p}\no
&+C\left(\|\na\uu\|_{L^p}+\|\uu\bar\xx^{-1/4}\|_{L^\infty} \|\bar \xx^{-3/2}\|_{L^p}+\|\na^2\uu\|_{L^p}\right)\|\ti\hh\|_{L^\infty}\no
\leq& C\left(\phi^\alpha(t)+\|\na^2\uu\|_{L^2\cap L^p}\right) \left(1+\|\na\ti\hh\|_{L^p}+\|\na\ti\hh\|_{L^q}\right).\nonumber
\end{align}
Combining \eqref{3.37} and \eqref{3.38}, one yields
\begin{align}\la{3.39}
&\frac{d}{dt}\left(\|\ti\hh\|_{L^2}+\|\na\ti\hh\|_{L^p}\right)\\
\leq & C\left(\phi^\alpha(t)+\|\na^2\uu\|_{L^2\cap L^p}\right) \left(1+\|\na\ti\hh\|_{L^p}+\|\ti\hh\|_{L^2}+ \|\na\ti\hh\|_{L^q}\right).\nonumber
\end{align}
Similarly, with the help of \eqref{3.34} and  denoting $w\triangleq\n\bar\xx^a$, one obtains for $p\in [2, q]$ that
\be\la{3.40}
\frac{d}{dt}\left(\|w\|_{L^2}+\|\na w\|_{L^p}\right)\leq C\left(\phi^\alpha(t)+\|\na^2\uu\|_{L^2\cap L^p}\right) \left(1+\|\na w\|_{L^p}+ \|w\|_{L^2}+\|\na w\|_{L^q}\right),
\ee
where $w$ satisfying the following equality
\be\la{a3.40}
w_t+\uu\cdot\na w-aw\uu\cdot\na\ln\bar\xx+w\div\uu=0.
\ee

Next, we claim that
\be\la{3.41}
\int_0^t\left(\|\na^2\uu\|_{L^2\cap L^q}^{(q+1)/q}+s\|\na^2\uu\|_{L^2\cap L^q}^2\right)ds\leq C\exp\left(C\int_0^t\phi^\alpha(s)ds\right),
\ee
which together with \eqref{3.39}, \eqref{3.40} and Gronwall inequality gives \eqref{3.33}.

Now, to finish the proof of Lemma \ref{l3.3}, it solely has to prove \eqref{3.41}. In fact, on the one hand, it follows from \eqref{3.4}, \eqref{3.21} and \eqref{3.24} that
\begin{align}\la{3.42}
&\int_0^t\left(\|\na^2\uu\|_{L^2}^{5/3}+s\|\na^2\uu\|_{L^2}^2\right)ds\\
\leq & C\int_0^t\left(\|\sqrt\n\uu_t\|_{L^2}^2+\phi^\alpha(s)\right)ds +C\exp\left(C\int_0^t\phi^\alpha(s)ds\right)\int_0^t\phi^\alpha(s)ds\no
\leq & C\exp\left(C\int_0^t\phi^\alpha(s)ds\right).\nonumber
\end{align}
On the other hand, choosing $p=q$ in \eqref{3.20} and using \eqref{3.9}, \eqref{3.10} and \eqref{3.12}, we have
\begin{align}\la{3.43}
\|\na^2\uu\|_{L^q}\leq &C\left(\|\n\uu_t\|_{L^p}+ \|\n\uu\cdot \na\uu\|_{L^p}+\|\na P(\n)\|_{L^p}+\||\hh||\na \hh|\|_{L^p}\right)\\
\leq & C \|\n\uu_t\|_{L^2}^{2(q-1)/(q^2-2)}\|\n \uu_t\|_{L^{q^2}}^{(q^2-2q)/(q^2-2)}+C\|\n\uu\|_{L^{2q}} \|\na\uu\|_{L^{2q}}+C\phi^\alpha(t)\no
\leq & C \phi^\alpha(t)\left(\|\sqrt\n\uu_t\|_{L^2} ^{2(q-1)/(q^2-2)}\|\na\uu_t \|_{L^2}^{(q^2-2q)/(q^2-2)} +\|\sqrt\n\uu_t\|_{L^2}\right)\no
&+C\phi^\alpha(t)\left(1+\|\na^2\uu\|_{L^2}^{(q-1)/q}\right).\nonumber
\end{align}
Then, combined \eqref{3.43} with \eqref{3.5}, \eqref{3.24} and \eqref{3.42} leads to
\begin{align}\la{3.44}
&\int_0^t\|\na^2\uu\|_{L^q}^{(q+1)/q}\\
\leq & C\int_0^t \phi^\alpha(s) s^{-(q+1)/2q}\left(s\|\sqrt\n\uu_t\|_{L^2}^2\right)^{\frac{q^2-1} {q(q^2-2)}} \left(s\|\na\uu_t\|_{L^2}^2\right)^{\frac{(q-2)(q+1)} {2(q^2-2)}}ds\no
&+C\int_0^t\|\sqrt\n\uu_t\|_{L^2}^2ds+C\exp\left(C\int_0^t \phi^\alpha(s)ds\right)\no
\leq & C\exp\left(C\int_0^t \phi^\alpha(s)ds\right)\left( 1+\int_0^t\left( \phi^\alpha(s) +s^{-\frac{q^3+q^2-2q-1}{q^3+q^2-2q}}+s \|\na\uu_t\|_{L^2}^2 \right)ds\right)\no
\leq &C\exp\left(C\int_0^t \phi^\alpha(s)ds\right).\nonumber
\end{align}
and that
\begin{align}\la{3.45}
\int_0^ts\|\na^2\uu\|_{L^q}^2ds\leq & C\int_0^t\phi^\alpha(s)\left(s\|\sqrt\n\uu_t\|_{L^2}^2\right) ^{\frac{2(q-1)} {q^2-2}} \left(s\|\na\uu_t\|_{L^2}^2 \right) ^{\frac{q^2-2q} {q^2-2}}ds\\
&+C\int_0^t\phi^\alpha(s)s\|\sqrt\n\uu_t\|_{L^2}^2ds+C\int_0^t \phi^\alpha(s) s\left(1+\|\na^2\uu\|_{L^2}^{2(q-1)/q}\right)ds\no
\leq &C\sup_{0\leq s\leq t}\left(s\|\sqrt\n\uu_t\|_{L^2}^2\right) ^{\frac{2(q-1)} {q^2-2}}\int_0^t\left(\phi^\alpha(s)+s\|\na\uu_t \|_{L^2}^2 \right)ds\no
&+C\sup_{0\leq s\leq t}\left(s\|\sqrt\n\uu_t \|_{L^2}^2\right) \int_0^t\phi^\alpha(s)ds+C\int_0^t\left(\phi^\alpha(s) +s\|\na^2\uu\|_{L^2}^2\right)ds\no
\leq &C\exp\left(C\int_0^t \phi^\alpha(s)ds\right).\nonumber
\end{align}
Therefore, \eqref{3.42}, \eqref{3.44} and \eqref{3.45} give \eqref{3.41} and we complete the proof of Lemma \ref{l3.3}. \hfill $\Box$

Now, we complete the proof of Proposition \ref{prop}, which is a direct consequence of Lemma \ref{l3.1}, \ref{l3.2} and \ref{l3.3}.

{\it Proof of Proposition \ref{prop}.} It follows from \eqref{3.4} and \eqref{3.33} that
$$
\phi(t)\leq \exp\left(C\exp\left(C\int_0^t \phi^\alpha(s)ds\right)\right).
$$
Standard arguments thus show that for $\ti M\triangleq e^{Ce}$ and $T_0 \triangleq \min\{T_1, (CM^\alpha)^{-1}\}$,
$$
\sup_{0\leq t\leq T_0}\phi(t)\leq \ti M,
$$
which together with \eqref{3.4}, \eqref{3.21} and \eqref{3.41} leads to \eqref{3.3}. Then the proof of Proposition \ref{prop} is finished. \hfill $\Box$

\section{ A priori estimates (II)}
In this section, in addition to $\mu, \lambda, \ga, q, a, \eta_0, N_0$ and $C_0$, the generic positive constant $C$ may also depend on $\delta_0$, $\|\na^2\uu_0\|_{L^2}$, $\|\na^2\n_0\|_{L^q}$, $\|\na^2P(\n_0)\|_{L^q}$, $\|\na^2\hh_0\|_{L^q}$, $\|\bar\xx^{\delta_0}\na^2\n_0\|_{L^2}$, $\|\bar\xx^{\delta_0} \na^2P(\n_0)\|_{L^2}$, $\|\bar\xx^{\delta_0}\na^2\hh_0\|_{L^2}$ and $\|\mathbf{g}\|_{L^2}$.

\bl\la{l4.1}
It holds that
\be\la{4.1}
\sup_{0\leq t\leq T_0}\left(\|\bar\xx^{\delta_0}\na^2\n\|_{L^2} +\|\bar\xx^{\delta_0}\na^2P(\n)\|_{L^2}+\|\bar\xx^{\delta_0} \na^2\hh\|_{L^2}\right)\leq C.
\ee
\el
\pf First, due to \eqref{1.13}, \eqref{2.1} and \eqref{2.2}$_2$, defining $$
\sqrt\n\uu_t(t=0, \xx)\triangleq -\mathbf{g}-\sqrt{\n_0}\uu_0\cdot\na\uu_0,
$$
integrating \eqref{3.32} over $(0, T_0)$ and using \eqref{3.3} and \eqref{3.4}, we have
\be\la{4.2}
\sup_{0\leq t\leq T_0}\|\sqrt\n\uu_t\|_{L^2}^2+\int_0^{T_0} \|\na\uu_t\|_{L^2}^2dt\leq C,
\ee
which together with \eqref{3.3} and \eqref{3.21} leads to
\begin{align}\la{4.3}
\|\na^2\uu\|_{L^2}\leq C\phi^{1/2}(t)\|\sqrt\n\uu_t\|_{L^2}+C\phi^\alpha (t), \nonumber\\
\sup_{0\leq t\leq T_0}\|\na\uu\|_{H^1}\leq C.
\end{align}
This combined with \eqref{3.3} and \eqref{3.35} shows that for $\delta\in (0, 1]$,
\be\la{4.4}
\|\n^\delta\uu\|_{L^\infty}+\|\bar\xx^{-\delta}\uu\|_{L^\infty}\leq C(\delta).
\ee
Directly calculate  that for $2\leq r\leq q$
\be\la{4.5}
\|\n_t(\bar\xx^{(1+a)/2}+|\uu|)\|_{L^r}+ \|P_t(\n)(1+|\uu|) \|_{L^r}+\|\hh_t(\bar\xx^{(1+a)/2} +|\uu|)\|_{L^r}\leq C,
\ee
due to \eqref{2.2}$_1$, \eqref{2.2}$_3$, \eqref{3.3}, \eqref{3.14}, \eqref{4.3} and \eqref{4.4}. It follows from \eqref{3.9}, \eqref{4.2}, \eqref{4.3} and \eqref{4.4} that for $\delta\in (0, 1]$ and $s>2/\delta$,
\begin{align}\la{4.6}
\|\bar\xx^{-\delta}\uu_t\|_{L^s}+\|\bar\xx^{-\delta}\uu\cdot \na\uu\|_{L^s}\leq & C\|\bar\xx^{-\delta}\uu_t\|_{L^s} +C\|\bar\xx^{-\delta}\uu\|_{L^\infty}\|\na\uu\|_{L^s}\\
\leq &C(\delta, s)+C(\delta, s)\|\na\uu_t\|_{L^2}.\nonumber
\end{align}

Next, denoting  $\bar\hh\triangleq \bar\xx^{\delta_0}\hh$ and $v \triangleq \bar\xx^{\delta_0}f(\n)$ with $f(\n)=\n^p$ for $p\in [1, \ga]$, we easily get from \eqref{3.3} that
\be\la{4.7}
\|v\|_{L^\infty}+\|\na v\|_{L^2\cap L^q}+\|\bar\hh\|_{L^\infty}+ \|\na \bar\hh\|_{L^2\cap L^q}\leq C,
\ee
where $v$ satisfies
\be\la{a4.7}
v_t+\uu\cdot\na v-\delta_0 v\uu\cdot\na \ln\bar\xx+p v\div\uu=0.
\ee
It follows from \eqref{2.2}$_3$, similarly as \eqref{3.36}, we have
\be\la{4.8}
\bar\hh_t+\uu\cdot\na\bar\hh-\delta_0\uu\cdot\na\ln\bar\xx\cdot \bar\hh +\bar\hh\div\uu=\bar\hh\cdot\na\uu.
\ee
Therefore, direct calculations give that
\begin{align}\la{4.9}
\frac{d}{dt}\|\na^2\bar\hh\|_{L^2}\leq & C\left(1+\|\na\uu\| _{L^\infty} +\|\uu\cdot\na\ln\bar\xx\|_{L^\infty}\right)\|\na^2\bar\hh\|_{L^2} +C\||\na^2\uu||\na\bar\hh|\|_{L^2}\\
&+C\||\na\bar\hh||\na\uu||\na\ln\bar\xx|\|_{L^2}+C\||\na\bar\hh||\uu| |\na^2 \ln\bar\xx|\|_{L^2}\no
&+C\|\bar\hh\|_{L^\infty}\left(\|\na^2\left(\uu\cdot\na\ln\bar\xx \right)\|_{L^2}+\|\na^3\uu\|_{L^2}\right)\no
\leq &C\left(1+\|\na\uu\|_{L^\infty}\right)\|\na^2\bar\hh\|_{L^2} +C\|\na^2 \uu\|_{L^{2q/(q-2)}}\|\na\bar\hh\|_{L^q} \no
&+C\|\na\bar\hh\|_{L^2}\|\na\uu\|_{L^\infty}+C\|\na\bar\hh\|_{L^2} \||\uu||\na^2\ln\bar\xx|\|_{L^\infty}\no
&+C\|\na^2\uu\|_{L^2}+C\|\na\uu\|_{L^2}+C\||\uu||\na^3\ln\bar\xx|\|_{L^2} +C\|\na^3\uu\|_{L^2}\no
\leq & C\left(1+\|\na\uu\|_{L^\infty}\right)\|\na^2\bar\hh\|_{L^2}+C +C\|\na^3\uu\|_{L^2},\nonumber
\end{align}
where in the second and third inequalities we have used \eqref{4.4} and \eqref{4.7}. Similarly, we can also obtain from \eqref{a4.7} after calculations
\be\la{4.10}
\frac{d}{dt}\|\na^2v\|_{L^2}\leq C\left(1+\|\na\uu\|_{L^\infty}\right)\|\na^2v\|_{L^2}+C +C\|\na^3\uu\|_{L^2}.
\ee
Combing \eqref{4.9} with \eqref{4.10}, we get
\begin{align}\la{4.12}
&\frac{d}{dt}\left(\|\na^2v\|_{L^2}+\|\na^2\bar\hh\|_{L^2}\right) \\
\leq& C\left(1+\|\na\uu\|_{L^\infty}\right)\left(\|\na^2v\|_{L^2}+ \|\na^2\bar\hh\|_{L^2}\right) +C+C\|\na^3\uu\|_{L^2}.\nonumber
\end{align}
We use \eqref{2.9} and \eqref{4.3} to estimate the last term on the right-hand side of \eqref{4.12} as follows:
\begin{align}\la{4.13}
\|\na^3\uu\|_{L^2}\leq &C \|\na(\n\uu_t)\|_{L^2}+C\|\na(\n\uu\cdot \na\uu)\|_{L^2}+C \|\na^2P(\n)\|_{L^2}\\
&+C\|\na^2|\hh|^2\|_{L^2}+C \|\na(\hh\cdot\na\hh)\|_{L^2}\no
\leq &C\|v\|_{L^\infty}\|\na\uu_t\|_{L^2}+C\|\bar\xx^a\na\n\|_{L^q}\|\bar \xx^{-a}\uu_t\|_{L^{2q/(q-2)}}\no
&+C\|\bar\xx^{-\delta}\uu\|_{L^\infty} \|\bar\xx^\delta\na\n\|_{L^q} \|\na\uu\|_{L^{2q/(q-2)}}+C\|v\|_{L^\infty}\|\na\uu\|_{L^4}^2\no
&+C\|v\|_{L^\infty}\|\bar\xx^{-\delta}\uu\|_{L^\infty}\|\na^2\uu\| _{L^2}^2 +C \|\na^2P(\n)\|_{L^2}\no
&+C\|\na\bar\hh\|_{L^4}^2 +C\|\bar\hh\|_{L^\infty}\|\na^2
\hh\|_{L^2}\no
\leq &C\|\na\uu_t\|_{L^2}+C\left(\|\na^2v\|_{L^2} +\|\na^2\bar\hh\|_{L^2}\right)+C,\nonumber
\end{align}
where in the last inequality we have used \eqref{3.3}, \eqref{4.3}, \eqref{4.4}, \eqref{4.6} and the following fact:
\begin{align}\la{4.14}
&\|\bar\xx^{\delta_0}\na^2\n\|_{L^2}+\|\bar\xx^{\delta_0}\na^2 P(\n)\|_{L^2}+\|\bar\xx^{\delta_0}\na^2 \hh\|_{L^2}\no
\leq&C\|\na^2(\bar\xx^{\delta_0}\n)\|_{L^2}+C\|\na^2 (\bar\xx^{\delta_0}P(\n))\|_{L^2}+\|\na^2(\bar\xx^{\delta_0} \hh)\|_{L^2}+C.
\end{align}
Substituting \eqref{4.13} into \eqref{4.12} and noting the definition of $v$, one has
\begin{align*}
&\frac{d}{dt}\left(\|\na^2(\bar\xx^{\delta_0}\n)\|_{L^2} +\|\na^2(\bar\xx^{\delta_0}P(\n))\|_{L^2} +\|\na^2(\bar\xx^{\delta_0}\hh)\|_{L^2}\right)\no
\leq &C(1+\|\na^2\uu\|_{L^q})\left(\|\na^2(\bar\xx^{\delta_0}\n)\|_{L^2} +\|\na^2(\bar\xx^{\delta_0}P(\n))\|_{L^2} +\|\na^2(\bar\xx^{\delta_0}\hh)\|_{L^2}\right)\no
&+C\|\na\uu_t\|_{L^2}+C,
\end{align*}
which together with \eqref{3.3}, \eqref{4.2}, \eqref{4.14} and Gronwall inequality gives \eqref{4.1} and completes the proof of Lemma \ref{l4.1}. \hfill $\Box$

\bl\la{l4.2}
It holds that
\be\la{4.15}
\sup_{0\leq t\leq T_0}t\|\na\uu_t\|_{L^2}^2+\int_0^{T_0}t\left( \|\sqrt\n \uu_{tt}\|_{L^2}^2+\|\na^2\uu_{t}\|_{L^2}^2\right)dt\leq C.
\ee
\el
\pf Multiplying \eqref{3.25} by $\uu_{tt}$ and integrating the resultant equality over $B_R$,  integration by parts lead to
\begin{align}\la{4.16}
&\frac12\frac{d}{dt}\left(\mu\|\na\uu_t\|_{L^2}^2+(\mu+\lambda) \|\div\uu_t\|_{L^2}^2\right)+\|\sqrt\n \uu_{tt}\|_{L^2}^2\\
=&-\int\left(2\n\uu\cdot\na\uu_t\cdot\uu_{tt}+\n\uu_t\cdot\na \uu \cdot\uu_{tt}\right)d\xx-\int\n\uu\cdot\na(\uu\cdot\na\uu)\cdot \uu_{tt}d\xx\no
&-\int\n\uu\cdot\na\uu_{tt}\cdot\uu_td\xx-\int\n\uu\cdot\na\uu_{tt} \cdot(\uu \cdot\na)\uu d\xx+\int P_t(\n)\div\uu_{tt}d\xx\no
&+\frac12\int |\hh|^2_t\div\uu_{tt}d\xx-\int\hh_t\cdot\na\uu_{tt}\cdot \hh d\xx-\int\hh\cdot\na\uu_{tt}\cdot\hh_td\xx.\nonumber
\end{align}
Now, we estimate each term on the right-hand side of \eqref{4.16}. First, it follows from \eqref{3.3}, \eqref{4.2}, \eqref{4.3}, \eqref{4.4} and \eqref{4.6} that
\begin{align}\la{4.17}
&\left|\int\left(2\n\uu\cdot\na\uu_t\cdot\uu_{tt}+\n\uu_t\cdot\na \uu \cdot\uu_{tt}\right)d\xx\right|+\left|\int\n\uu\cdot\na(\uu\cdot\na\uu) \cdot \uu_{tt}d\xx\right|\\
\leq &\ve\|\sqrt\n \uu_{tt}\|_{L^2}^2 +C(\ve)\left(\|\sqrt\n\uu\| _{L^\infty}^2 \|\na\uu_t\|_{L^2}^2+\|\sqrt\n\uu_t\|^2_{L^4}\|\na\uu\| _{L^4}^2\right)\no
&+C(\ve)\left(\|\sqrt\n\uu\|_{L^\infty}^2\|\na\uu\|_{L^4}^2 +\|\n^{1/4}\uu\|_{L^\infty}^4\|\na^2\uu\|_{L^2}^2\right) \no
\leq &\ve\|\sqrt\n \uu_{tt}\|_{L^2}^2 +C(\ve)\left(1+\|\na\uu_t\|_{L^2}^2 \right).\nonumber
\end{align}

Next, direct calculations yield that
\begin{align}\la{4.18}
&-\int\n\uu\cdot\na\uu_{tt}\cdot\uu_td\xx-\int\n\uu\cdot\na\uu_{tt} \cdot(\uu \cdot\na)\uu d\xx\\
=&-\frac{d}{dt}\int\left( \n\uu\cdot\na\uu_{t}\cdot\uu_t+\n\uu \cdot\na\uu_{t}\cdot(\uu \cdot\na)\uu\right)d\xx +\int(\n\uu)_t\cdot\na\uu_{t}\cdot\uu_td\xx\no
&+\int(\n\uu)_t\cdot\na\uu_{t}\cdot(\uu \cdot\na)\uu d\xx+\int\n\uu\cdot\na\uu_{t}\cdot\uu_{tt}d\xx\no
&+\int\n\uu\cdot\na\uu_{t}\cdot(\uu_t \cdot\na)\uu d\xx +\int\n\uu\cdot\na\uu_{t}\cdot(\uu \cdot\na)\uu_t d\xx.\nonumber
\end{align}
First, H\"{o}lder inequality together with \eqref{3.3} and \eqref{4.4}--\eqref{4.6} gives
\begin{align}\la{4.19}
\int(\n\uu)_t\cdot\na\uu_{t}\cdot\uu_td\xx
\leq & C\|\n\bar\xx^a\|_{L^\infty}\|\bar\xx^{-a/2}\uu_t\|_{L^4}^2 \|\na\uu_t\|_{L^2}\\
&+C\|\bar\xx^{(1+a)/2}\n_t\|_{L^2} \|\bar\xx^{-1/2}\uu\|_{L^\infty}\|\bar\xx^{-a/2}\uu_t\|_{L^4} \|\na\uu_t\|_{L^4}\no
\leq &\delta\|\na^2\uu_t\|_{L^2}^2+C(\delta)\|\na\uu_t\|_{L^2}^4 +C(\delta).\nonumber
\end{align}
Similarly,
\begin{align}\la{4.20}
&\int(\n\uu)_t\cdot\na\uu_{t}\cdot(\uu \cdot\na)\uu d\xx\\
\leq  & C\|\n\bar\xx^a\|_{L^\infty} \|\bar\xx^{-a/2}\uu_t\|_{L^4} \|\bar\xx^{-a/2}\uu\cdot\na\uu\|_{L^4}\|\na\uu_t\|_{L^2}\no
&+C\|\bar\xx^{(1+a)/2}\n_t\|_{L^2}\|\bar\xx^{-1/2}\uu\|_{L^\infty} \|\bar\xx^{-a/2}\uu\|_{L^\infty}\|\na\uu\|_{L^4}\|\na\uu_t\|_{L^4}\no
\leq &\delta\|\na^2\uu_t\|_{L^2}^2+C(\delta)\|\na\uu_t\|_{L^2}^4 +C(\delta).\nonumber
\end{align}
Then, H\"{o}lder inequality together with \eqref{4.4} leads to
\begin{align}\la{4.21}
\int\n\uu\cdot\na\uu_{t}\cdot\uu_{tt}d\xx\leq & C\|\sqrt\n \uu_{tt}\|_{L^2}\|\sqrt\n\uu\|_{L^\infty}\|\na\uu_t\|_{L^2}\\
\leq &\ve\|\sqrt\n \uu_{tt}\|_{L^2}^2 +C(\ve)\|\na\uu_t\|_{L^2}^2.\nonumber
\end{align}
Next, it follows from \eqref{4.3}, \eqref{4.4} and \eqref{4.6} that
\begin{align}\la{4.22}
\int\n\uu\cdot\na\uu_{t}\cdot(\uu_t \cdot\na)\uu d\xx\leq &C\|\n\bar\xx^a\|_{L^\infty}\|\bar\xx^{-a/2} \uu\|_{L^\infty} \|\bar\xx^{-a/2} \uu_t\|_{L^4}\|\na\uu\|_{L^4}\|\na\uu_t\|_{L^2}\no
\leq & C+C(\ve)\|\na\uu_t\|_{L^2}^2,
\end{align}
and similarly,
\begin{align}\la{4.23}
\int\n\uu\cdot\na\uu_{t}\cdot(\uu \cdot\na)\uu_t d\xx \leq &C\|\n\bar\xx^a\|_{L^\infty}\|\bar\xx^{-a/2} \uu\|_{L^\infty}^2\|\na\uu_t\|_{L^2}^2\\
\leq & C\|\na\uu_t\|_{L^2}^2.\nonumber
\end{align}
Inserting \eqref{4.19}--\eqref{4.23} into \eqref{4.18} shows
\begin{align}\la{4.24}
&-\int\n\uu\cdot\na\uu_{tt}\cdot\uu_td\xx-\int\n\uu\cdot\na\uu_{tt} \cdot(\uu \cdot\na)\uu d\xx\\
=&-\frac{d}{dt}\int\left( \n\uu\cdot\na\uu_{t}\cdot\uu_t+\n\uu \cdot\na\uu_{t}\cdot(\uu \cdot\na)\uu\right)d\xx+\ve\|\sqrt\n \uu_{tt}\|_{L^2}^2\no
&+C(\ve, \delta)\|\na\uu_t\|_{L^2}^4+\delta\|\na^2\uu_t\|_{L^2}^2 +C(\ve, \delta).\nonumber
\end{align}

Next, it follows from \eqref{3.14}, \eqref{4.3} and \eqref{4.5} that
\begin{align}\la{4.25}
\int P_t(\n)\div\uu_{tt}d\xx=&\frac{d}{dt}\int P_t(\n)\div\uu_t d\xx-\int (P(\n)\uu)_t\cdot\na\div\uu_t d\xx\\
&+(\ga-1)\int(P(\n)\div\uu)_t\div\uu_t d\xx\no
\leq &\frac{d}{dt}\int P_t(\n)\div\uu_t d\xx +C\left(\|P_t(\n)\uu\|_{L^2}+\|P(\n)\uu_t\|_{L^2}\right) \|\na^2\uu_t\|_{L^2}\no
&+C\|P_t(\n)\|_{L^2}\|\na\uu\|_{L^4}\|\na\uu_t\|_{L^4} +\|P(\n)\|_{L^\infty}\|\na\uu_t\|_{L^2}^2\no
\leq &\frac{d}{dt}\int P_t(\n)\div\uu_t d\xx +\delta\|\na^2 \uu_t\|_{L^2}^2+C(\delta)\left(\|\na\uu_t\|_{L^2}^2+1\right).\nonumber
\end{align}
Similarly, \eqref{2.2}$_3$, \eqref{4.1}, \eqref{4.3} and \eqref{4.5} give that
\begin{align}\la{4.26}
&\frac12\int |\hh|^2_t\div\uu_{tt}d\xx-\int\hh_t\cdot\na\uu_{tt}\cdot \hh d\xx-\int\hh\cdot\na\uu_{tt}\cdot\hh_td\xx\\
=&\frac{d}{dt}\int\left(\frac12|\hh|^2_t\div\uu_t-\hh_t\cdot\na\uu _t\cdot\hh -\hh\cdot\na\uu_t\cdot\hh_t\right)d\xx-\int|\hh_t|^2 \div\uu_td\xx \no
&-\int\hh\cdot\hh_{tt}\div\uu_td\xx+\int\hh_{tt}\cdot\na\uu_t\cdot\hh d\xx+2\int\hh_t\cdot\na\uu_t\cdot\hh_t d\xx\no
&+\int\hh\cdot\na\uu_t\cdot \hh_{tt}d\xx\no
\leq &\frac{d}{dt}\int\left(\frac12|\hh|^2_t\div\uu_t-\hh_t\cdot\na\uu _t\cdot\hh -\hh\cdot\na\uu_t\cdot\hh_t\right)d\xx\no
&+C\|\hh_t\bar\xx^{(a+1)/2}\|_{L^q} \|\na\uu_t\|_{L^2}\left(\|\uu\bar\xx^{-1/2}\|_{L^\infty} \|\na\hh\|_{(q-2)/2q}+\|\hh\|_{L^\infty}\|\na\uu\|_{(q-2)/2q}\right)\no
&+C\|\bar\xx^{a/2}\hh\|_{L^\infty}\|\bar\xx^{-a/2}\uu_t\|_{L^4} \|\na\hh\|_{L^4} \|\na\uu_t\|_{L^2}+C\|\hh\|_{L^\infty}^2\|\na \uu_t\|_{L^2}^2\no
&+C\|\bar\xx^{a/2}\hh\|_{L^\infty} \|\bar\xx^{-a/2}\uu\|_{L^\infty} \|\na\hh_t\|_{L^2} \|\na\uu_t\|_{L^2}\no
&+C\|\hh\|_{L^\infty}\|\hh_t \bar\xx^{a/2}\|_{L^q}\|\na\uu_t\|_{L^2}\|\na\uu\|_{(q-2)/2q}\no
\leq &\frac{d}{dt}\int\left(\frac12|\hh|^2_t\div\uu_t-\hh_t\cdot\na\uu _t\cdot\hh -\hh\cdot\na\uu_t\cdot\hh_t\right)d\xx+C \left(\|\na\uu_t\|_{L^2}^2+1\right),\nonumber
\end{align}
where in the last inequality we have used \eqref{4.6} and the following simple fact:
\begin{align}\la{4.27}
\|\na\hh_t\|_{L^2}\leq& C\||\na\uu||\na\hh|\|_{L^2}+C\||\uu||\na^2\hh|\| _{L^2}+C\||\hh||\na^2\uu|\|_{L^2}\\
\leq &C\|\na\uu\|_{L^{2q/(q-2)}}\|\na\hh\|_{L^q}+C\|\bar\xx^{-\delta_0} \uu\|_{L^\infty} \|\bar\xx^{\delta_0}\na^2\hh\|_{L^2}+C\|\hh\|_{L^\infty} \|\na^2\uu\|_{L^2}\no
\leq &C,\nonumber
\end{align}
due to \eqref{4.1}, \eqref{4.3}, \eqref{4.4} and \eqref{4.7}.

Inserting \eqref{4.17} and \eqref{4.24}--\eqref{4.26} into \eqref{4.16} and choosing $\ve$ suitably small we follow that
\be\la{4.28}
\Phi'(t)+\|\sqrt\n \uu_{tt}\|_{L^2}^2\leq C\delta\|\na^2\uu_t\|_{L^2}^2 +C\|\na\uu_t\|_{L^2}^4+C,
\ee
where
\begin{align*}
\Phi(t)\triangleq&\mu\|\na\uu_t\|_{L^2}^2+(\mu+\lambda) \|\div\uu_t\|_{L^2}^2+\int\left( \n\uu\cdot\na\uu_{t}\cdot\uu_t+\n\uu \cdot\na\uu_{t}\cdot(\uu \cdot\na)\uu\right)d\xx\no
&-\int P_t(\n)\div\uu_t d\xx-\int\left(\frac12|\hh|^2_t\div\uu_t-\hh_t\cdot\na\uu _t\cdot\hh -\hh\cdot\na\uu_t\cdot\hh_t\right)d\xx
\end{align*}
satisfies
\be\la{4.29}
C(\mu)\|\na\uu_t\|_{L^2}^2-C\leq \Phi(t)\leq C\|\na\uu_t\|_{L^2}^2 +C,
\ee
due to the following fact:
\begin{align*}
&\left|\int\left( \n\uu\cdot\na\uu_{t}\cdot\uu_t+\n\uu \cdot\na\uu_{t}\cdot(\uu \cdot\na)\uu\right)d\xx-\int P_t(\n)\div\uu_t d\xx\right|\no
&+\left|\int\left(\frac12|\hh|^2_t\div\uu_t-\hh_t\cdot\na\uu _t\cdot\hh -\hh\cdot\na\uu_t\cdot\hh_t\right)d\xx\right|\no
\leq &C\|\sqrt\n \uu\|_{L^\infty}\|\na\uu_t\|_{L^2}\|\sqrt\n\uu_t\|_{L^2}
+C\|\sqrt\n \uu\|_{L^\infty}^2\|\na\uu_t\|_{L^2}\|\na\uu\|_{L^2}\no
&+C\|\hh\|_{L^\infty}\|\hh_t\|_{L^2}\|\na\uu_t\|_{L^2}\no
\leq &\ve\|\na\uu_t\|_{L^2}^2+C(\ve),
\end{align*}
which yields from \eqref{4.2}--\eqref{4.5}.

Then, it remains to estimate the first term on the right-hand side of \eqref{4.27}. In fact, it follows from \eqref{3.25} and Lemma \ref{l2.5} that, for $s>2$,
\begin{align}\la{4.30}
&\|\na^2\uu_t\|_{L^2}^2\\
\leq & C\|\n\|_{L^\infty} \|\sqrt\n\uu_{tt} \|_{L^2}^2 +C\|\n\| _{L^\infty} \|\sqrt\n\uu \|_{L^\infty}^2 \|\na\uu_t\|_{L^2}^2 +C\|\bar\xx^{(a+1)/2}\n_t\|_{L^q}^2\|\bar\xx^{-1} \uu_t\|_{L^{2q/(q-2)}}^2 \no &+C\|\bar\xx^{(a+1)/2}\n_t\|_{L^q}^2\|\bar\xx^{-1}\uu\|_{L^\infty}^2  \|\na\uu\|_{L^{2q/(q-2)}}^2+C\|\bar\xx^{(a+1)/2}\n\|_{L^q}^2 \|\bar\xx^{-1}\uu_t\|_{L^s}^2\|\na\uu\|_{L^{2qs/(qs-2(q+s))}}^2\no
&+C\|\na P_t(\n)\|_{L^2}^2+C\|\bar\xx^{(a+1)/2}\hh_t\|_{L^q}^2 \|\na\hh\|_{L^{2q/(q-2)}}^2+C\|\hh\|_{L^\infty}^2\|\na\hh_t\|_{L^2}^2\no
\leq & C\|\sqrt\n\uu_{tt} \|_{L^2}^2+C\|\na\uu_t\|_{L^2}^4+C,\nonumber
\end{align}
due to \eqref{4.1}--\eqref{4.7}, \eqref{4.27} and the following fact:
\begin{align*}
\|\na P_t(\n)\|_{L^2}\leq &C\|\na\uu\|_{L^{2q/(q-2)}}\|\na v\|_{L^q} +C\|\bar\xx^{-\delta_0}\uu\|_{L^\infty}\|\bar\xx^{\delta_0}\na^2 P(\n)\|_{L^2} +C\|\na^2\uu\|_{L^2}\\
\leq & C,
\end{align*}
where we have used Gagliardo-Nirenberg inequality, \eqref{4.1},  \eqref{4.3}, \eqref{4.4} and \eqref{4.7}.
Inserting \eqref{4.30} into \eqref{4.28} and choosing $\delta$ suitably small  lead to
\be\la{4.31}
\Phi'(t)+\|\sqrt\n \uu_{tt}\|_{L^2}^2\leq C\|\na\uu_t\|_{L^2}^4+C.
\ee

Multiplying \eqref{4.31} by $t$ and integrating the resultant inequality over $(0, T_0)$, we obtain from Gronwall inequality, \eqref{4.2} and \eqref{4.29}
$$
\sup_{0\leq t\leq T_0}t\|\na\uu_t\|_{L^2}^2+\int_0^{T_0}t\|\sqrt\n \uu_{tt}\|_{L^2}^2dt\leq C,
$$
which together with \eqref{4.30} yields \eqref{4.15} and completes the proof of Lemma \ref{l4.2}. \hfill $\Box$

\bl\la{l4.3}
It holds that
\be\la{4.32}
\sup_{0\leq t\leq T_0}\left( \|\na^2\n\|_{L^q}+\|\na^2P(\n)\|_{L^q} +\|\na^2\hh\|_{L^q}\right)\leq C.
\ee
\el
\pf Applying the differential operator $\na^2$ to \eqref{a4.7} and \eqref{4.8}, respectively, and multiplying each equality by $q|\na^2v| ^{q-2} \na^2v$ and $q|\na^2\bar\hh|^{q-2}\na^2\bar\hh$, and integrating the resultant equalities over $B_R$ lead to
\begin{align}\la{4.33}
&\frac{d}{dt}\left(\|\na^2\bar\hh\|_{L^q}+\|\na^2v\|_{L^q}\right)\\
\leq &C\|\na\uu\|_{L^\infty}\left(\|\na^2\bar\hh\|_{L^q} +\|\na^2 v\|_{L^q}\right)+ \left( \|\na\bar\hh\| _{L^\infty} +\|\na v\|_{L^\infty}\right)\|\na^2\uu\|_{L^q}\no
\leq &C\left(1+\|\na^2\uu\|_{L^q}\right)\left(1+\|\na^2\bar\hh\|_{L^q} +\|\na^2 v\|_{L^q}\right)+C\|\na^3\uu \|_{L^q}.\nonumber
\end{align}
Due to \eqref{2.9}, the last term on the right-hand side of \eqref{4.33} can be estimated as follows:
\begin{align}\la{4.34}
\|\na^3\uu\|_{L^q}\leq &C\|\na(\n\uu_t)\|_{L^q}+C\|\na(\n\uu\cdot\na \uu)\|_{L^q}+C\|\na^2 P(\n)\|_{L^q}\\
&+C\|\na ((\na\times\hh)\times\hh )\|_{L^q}\no
\leq & C\|\bar\xx^{-a}\uu_t\|_{L^\infty}\|\bar\xx^a\na\n\|_{L^q} +C\|\bar\xx^a\n\|_{L^\infty}\|\bar\xx^{-a}\na\uu_t\|_{L^q}\no
&+C\|\bar\xx^a\na\n\|_{L^q}\|\bar\xx^{-a}\uu\|_{L^\infty} \|\na\uu\|_{L^\infty}+C\|\n\|_{L^\infty}\|\na\uu\|_{L^{2q}}^2 \no &+C\|\bar\xx^a\n\|_{L^\infty}\|\bar\xx^{-a}\uu\|_{L^\infty} \|\na^2\uu\|_{L^q}+C\|\na^2 P(\n)\|_{L^q}\no
&+C\|\bar\xx^{a}\hh\|_{L^\infty}\|\bar\xx^{-a}\na^2\hh\|_{L^q} +C\|\na\hh\|_{L^{2q}}^2\no
\leq &C\|\na\uu_t\|_{L^q}+C\|\bar\xx^{-a}\uu_t\|_{L^q} +\frac12\|\na^3\uu\|_{L^q}+C\|\na^2 P(\n)\|_{L^q}\no
&+C\|\na^2\bar\hh\|_{L^q}+C\no
\leq &C\|\na\uu_t\|_{L^2}^{2/q}\|\na^2\uu_t\|_{L^2}^{(q-2)/q} +C\|\na\uu_t\|_{L^2}+\frac12\|\na^3\uu\|_{L^q}\no
&+C\|\na^2 P(\n)\|_{L^q}+C\|\na^2\bar\hh\|_{L^q}+C,\nonumber
\end{align}
where we have used \eqref{3.3}, \eqref{4.3}, \eqref{4.4}, \eqref{4.6} and \eqref{4.7}.

Next, it follows from \eqref{4.15} that
\begin{align}\la{4.35}
&\int_0^{T_0}\left(\|\na\uu_t\|_{L^2}^{2/q}\|\na^2\uu_t\|_{L^2} ^{(q-2)/q}\right)^{(q+1)/q}dt\\
\leq &C\sup_{0\leq t\leq T_0}\left(t\|\na\uu_t\|_{L^2}^2\right) ^{(q+1)/q^2}\int_0^{T_0}\left(t\|\na^2\uu_t\|_{L^2}^2 +t^{-(q^2+q)/(q^2+q+2)}\right)dt\no
\leq &C.\nonumber
\end{align}
Putting \eqref{4.34} into \eqref{4.33}, we get \eqref{4.32} from Gronwall inequality, \eqref{3.3}, \eqref{4.2} and \eqref{4.35}. Therefore, the proof of Lemma \ref{l4.3} is finished. \hfill $\Box$

\bl\la{l4.4}
It holds that
\begin{align}\la{4.36}
&\sup_{0\leq t\leq T_0}t\left( \|\na^3\uu\|_{L^2\cap L^q}+\|\na\uu_t\|_{H^1} +\|\na^2(\n\uu)\|_{L^{(q+2)/2}}\right)\\
&\qquad\qquad\quad+\int_0^{T_0}t^2\left(\|\na\uu_{tt}\| _{L^2}^2+\|\bar\xx^{-1}\uu_{tt}\|_{L^2}^2 \right)dt \leq C.\nonumber
\end{align}
\el
\pf First, we claim that
\be\la{4.37}
\sup_{0\leq t\leq T_0}t^2\|\sqrt\n\uu_{tt}\|_{L^2}^2+\int_0^{T_0} t^2\|\na\uu_{tt}\|_{L^2}^2 dt\leq C,
\ee
which together with \eqref{2.5}, \eqref{4.15} and \eqref{4.30} yields that
\be\la{4.38}
\sup_{0\leq t\leq T_0}t\|\na\uu_t\|_{H^1}+\int_0^{T_0}t^2\|\bar\xx^{-1} \uu_{tt}\|_{L^2}^2dt \leq C.
\ee
This combined with \eqref{4.13}, \eqref{4.32}, \eqref{4.34} and \eqref{4.35} leads to
\be\la{4.39}
\sup_{0\leq t\leq T_0}t\|\na^3\uu\|_{L^2\cap L^q}\leq C,
\ee
which together with \eqref{3.3}, \eqref{4.1} and \eqref{4.32}, shows
\begin{align}\la{4.40}
t\|\na^2(\n\uu)\|_{L^{(q+2)/2}}\leq &C t\||\na^2\n||\uu|\|_{L^{(q+2)/2}} +Ct\||\na\n||\na\uu|\|_{L^{(q+2)/2}}+Ct\|\n|\na^2\uu|\|_{L^{(q+2)/2}}\no
\leq &Ct\|\bar\xx^{\delta_0}\na^2\n\|_{L^2}^{2/(q+2)}\|\na^2\n\|_{L^q} ^{q/(q+2)}\|\bar\xx^{-2\delta_0/(q+2)}\uu\|_{L^\infty}\\
&+Ct\|\na\n\|_{L^q}\|\na\uu\|_{L^{q(q+2)/(q-2)}}+Ct\|\na^2\uu\| _{L^{(q+2)/2}}\no
\leq &C.\nonumber
\end{align}
Therefore, we complete the proof of \eqref{4.36} from \eqref{4.37}--\eqref{4.40}.

Now, we focus on the estimates of \eqref{4.37}. In fact, differentiating \eqref{3.25} with respect to $t$ yields that
\begin{align*}
&\n\uu_{ttt}+\n\uu\cdot\na\uu_{tt}-\mu\triangle\uu_{tt}-(\mu+\lambda) \na\div \uu_{tt}\no
=&2\div(\n\uu)\uu_{tt}+\div(\n\uu)_t\uu_t-2(\n\uu)_t\cdot\na\uu_t -\n_{tt}\uu\cdot\na\uu-2\n_t\uu_t\cdot\na\uu\no
&-\n\uu_{tt}\cdot\na\uu-\na P_{tt}(\n)-\frac12\na |\hh|_{tt}^2+\hh_{tt}\cdot\na\hh+2\hh_t\cdot\na\hh_t+\hh\cdot\na\hh_{tt},
\end{align*}
which multiplied by $\uu_{tt}$ and integrated by parts over $B_R$, shows that
\begin{align}\la{4.41}
&\frac12\frac{d}{dt}\int\n|\uu_{tt}|^2d\xx+\int\left(\mu|\na\uu_{tt}|^2 +(\mu+\lambda)|\div\uu_{tt}|^2\right)d\xx\\
=&-4\int\n\uu\cdot\na\uu_{tt}\cdot\uu_{tt}d\xx-\int(\n\uu)_t\cdot\left( \na(\uu_t\cdot\uu_{tt})+2\na\uu_t\cdot\uu_{tt}\right)\no
&-\int(\n\uu)_t\cdot\na(\uu\cdot\na\uu\cdot\uu_{tt})d\xx-2\int\n_t\uu_t \cdot \na\uu\cdot\uu_{tt}d\xx\no
&-\int\n\uu_{tt}\cdot\na\uu\cdot\uu_{tt}d\xx+\int P_{tt}(\n)\div\uu_{tt} d\xx+\frac12\int |\hh|^2_{tt}\div\uu_{tt}d\xx\no
&+\int\hh_{tt}\cdot\na\hh\cdot\uu_{tt}d\xx+2\int\hh_t\cdot\na\hh_t\cdot \uu_{tt}d\xx+\int\hh\cdot\na\hh_{tt}\cdot\uu_{tt}d\xx\no
\triangleq&\sum_{i=1}^{10}I_i.\nonumber
\end{align}

Now, we'll estimate each term on the right-hand side of \eqref{4.41} as follows:

First, it follows from \eqref{4.4} that
\be\la{4.42}
|I_1|\leq C\|\sqrt\n\uu\|_{L^\infty}\|\sqrt\n\uu_{tt}\|_{L^2}\|\na \uu_{tt}\|_{L^2}\leq \ve\|\na \uu_{tt}\|_{L^2}^2+C(\ve) \|\sqrt\n\uu_{tt}\|_{L^2}^2.
\ee

Next, H\"{o}lder inequality leads to
\begin{align}\la{4.43}
|I_2|\leq &C\|\bar\xx(\n\uu)_t\|_{L^q}\left(\|\bar\xx^{-1}\uu_{tt}\| _{2q/(q-2)}\|\na\uu_t\|_{L^2}+\|\bar\xx^{-1}\uu_{t}\| _{2q/(q-2)}\|\na\uu_{tt}\|_{L^2}\right)\\
\leq &C\left(1+\|\na\uu_t\|_{L^2}^2\right)\left(\|\sqrt\n\uu_{tt}\| _{L^2}+\|\na\uu_{tt}\|_{L^2} \right)\no
\leq &\ve\left(\|\sqrt\n\uu_{tt}\| _{L^2}^2+\|\na\uu_{tt}\|_{L^2}^2 \right)+C(\ve)\left(1+\|\na\uu_t\|_{L^2}^4\right),\nonumber
\end{align}
where we have used \eqref{2.6}, \eqref{3.9} and \eqref{4.5} and the following facts:
\begin{align}\la{4.44}
\|\bar\xx(\n\uu)_t\|_{L^q}\leq &C\|\bar\xx|\n_t||\uu|\|_{L^q} +C\|\bar\xx \n|\uu_t|\|_{L^q}\\
\leq &C\|\n_t\bar\xx^{(1+a)/2}\|_{L^q}\|\bar\xx^{-(a-1)/2}\uu\| _{L^\infty} +C\|\n\bar\xx^a\|_{L^{2q/(3-\ti a)}}\|\uu_t\bar\xx^{1-a}\|_{L^{2q/(\ti a-1)}}\no
\leq &C+C\|\na\uu_t\|_{L^2}.\nonumber
\end{align}
by using \eqref{4.4}--\eqref{4.6}, where $\ti a=\min\{2, a\}$.

Then, it follows from \eqref{3.9}, \eqref{4.3}, \eqref{4.4} and \eqref{4.44} that
\begin{align}\la{4.45}
|I_3|\leq &C\int|(\n\uu)_t|\left(|\uu||\na^2\uu||\uu_{tt}|+|\uu||\na\uu| |\na\uu_{tt}|+|\na\uu|^2|\uu_{tt}|\right)d\xx\\
\leq &C\|\bar\xx(\n\uu)_t\|_{L^q}\|\bar\xx^{-1/q}\uu\|_{L^\infty}\|\na^2\uu\| _{L^2} \|\bar\xx^{-(q-1)/q}\uu_{tt}\|_{L^{2q/(q-2)}}\no
&+C\|\bar\xx(\n\uu)_t\|_{L^q}\|\bar\xx^{-1}\uu\|_{L^\infty}\|\na\uu\| _{L^{2q/ (q-2)}}\|\na\uu_{tt}\|_{L^2}\no
&+C\|\bar\xx(\n\uu)_t\|_{L^q}\|\na\uu\|_{L^4}^2\|\bar\xx^{-1}\uu_{tt}\| _{L^{2q/(q-2)}}\no
\leq &C(1+\|\na\uu_t\|_{L^2})\left(\|\sqrt\n\uu_{tt}\| _{L^2}+\|\na\uu_{tt}\|_{L^2} \right)\no
\leq &\ve\left(\|\sqrt\n\uu_{tt}\| _{L^2}^2+\|\na\uu_{tt}\|_{L^2}^2 \right)+C(\ve)\left(1+\|\na\uu_t\|_{L^2}^2\right).\nonumber
\end{align}

Next, it follows from Cauchy inequality, together with \eqref{3.9}, \eqref{4.5} and \eqref{4.6} that
\begin{align}\la{4.46}
|I_4|\leq &C\int|\n_t||\uu_t||\na\uu||\uu_{tt}|d\xx\\
\leq & C\|\bar\xx\n_t\|_{L^q}\|\bar\xx^{-1/2}\uu_t\|_{L^{4q/(q-2)}} \|\na\uu\|_{L^2}\|\bar\xx^{-1/2}\uu_{tt}\|_{L^{4q/(q-2)}}\no
\leq & C\left(1+\|\na\uu_t\|_{L^2}\right)\left(\|\sqrt\n\uu_{tt}\| _{L^2}+\|\na\uu_{tt}\|_{L^2} \right)\no
\leq &\ve\left(\|\sqrt\n\uu_{tt}\| _{L^2}^2+\|\na\uu_{tt}\|_{L^2}^2 \right)+C(\ve)\left(1+\|\na\uu_t\|_{L^2}^2\right).\nonumber
\end{align}

Then, Gagliardo-Nirenberg inequality together with \eqref{4.3} gives
\begin{align}\la{4.47}
|I_5|\leq C\|\na\uu\|_{L^\infty}\|\sqrt\n\uu_{tt}\|_{L^2}^2\leq C\left(1+\|\na^2\uu\|_{L^q}\right)\|\sqrt\n\uu_{tt}\|_{L^2}^2.
\end{align}

Next, it follows from \eqref{3.3}, \eqref{3.14} and \eqref{4.3}--\eqref{4.6} that
\begin{align}\la{a4.47}
&\|P_{tt}(\n)\|_{L^2}\\
\leq & C\||\uu_t||\na P(\n)|\|_{L^2} +C\||\uu||\na P_t(\n)|\|_{L^2}+C\||P_t(\n)| |\div\uu|\|_{L^2} +C\|P(\n)|\div\uu_t|\|_{L^2}\no
\leq & C\|\bar\xx^{-\theta_0}\uu_t\|_{2q/[(q-2)\theta_0]} \|\bar\xx^{\theta_0} \na P(\n)\|_{L^{2q/[q-(q-2)\theta_0]}}\no
&+C\|\bar\xx^{-\delta_0/2}\uu\|_{L^\infty}\|\bar\xx^{\delta_0/2}\na P_t(\n)\|_{L^2}+C\|P_t(\n)\|_{L^q}\|\na\uu\|_{L^{2q/(q-2)}}+C\|\na \uu_t\|_{L^2}\no
\leq &C\left(1+\|\na\uu_t\|_{L^2}\right),\nonumber
\end{align}
where in the last inequality we have used the following simple fact that
\begin{align*}
\|\bar\xx^{\delta_0/2}\na P_t(\n)\|_{L^2}\leq &C \|\bar\xx^{\delta_0/2} |\uu||\na^2 P(\n)|\|_{L^2}+C\|\bar\xx^{\theta_0}|\na\uu||\na P(\n)| \|_{L^2}\\
&+C\|\bar\xx^{a\ga} P(\n)\na^2\uu\|_{L^2}\\
\leq &C\|\bar\xx^{-\delta_0/2}\uu\|_{L^\infty}\|\bar\xx^{\delta_0}\na^2 P(\n)\|_{L^2}+C\|\na\uu\|_{L^{2q/(q-2)}}\|\bar\xx^{\theta_0}\na P(\n)\|_{L^q}\\
&+C\|\bar\xx^{a\ga}P(\n)\|_{L^\infty}\|\na^2\uu\|_{L^2}\\
\leq&C,
\end{align*}
due to \eqref{3.3}, \eqref{4.1}, \eqref{4.3} and \eqref{4.4}. Then, \eqref{a4.47} together with Cauchy inequality leads to
\begin{align}\la{4.48}
|I_6|\leq &C\|P_{tt}(\n)\|_{L^2}\|\div\uu_{tt}\|_{L^2}\leq \ve\|\na\uu_{tt}\|_{L^2}^2+C(\ve)\left(1+\|\na\uu_t\|_{L^2}^2\right).
\end{align}

Finally, it follows from \eqref{2.2}$_3$, \eqref{3.3} and \eqref{4.3}--\eqref{4.6} that
\begin{align}\la{4.49}
\|\hh_{tt}\|_{L^2}\leq & C\|\bar\xx^{-\theta_0}\uu_t\| _{2q/[(q-2)\theta_0]} \|\bar\xx^{\theta_0} \na \hh\|_{L^{2q/[q-(q-2)\theta_0]}}+C\|\na \uu_t\|_{L^2}\\
&+C\|\bar\xx^{-\delta_0/2}\uu\|_{L^\infty}\|\bar\xx^{\delta_0/2}\na \hh_t\|_{L^2}+C\|\hh_t\|_{L^q}\|\na\uu\|_{L^{2q/(q-2)}}\no
\leq &C\left(1+\|\na\uu_t\|_{L^2}\right),\nonumber
\end{align}
where in the last inequality we have used the following fact:
\begin{align*}
\|\bar\xx^{\delta_0/2}\na\hh_t\|_{L^2}\leq& C\|\bar\xx^{\delta_0/2}|\na\uu||\na\hh|\|_{L^2}+C\|\bar\xx^{\delta_0/2} |\uu||\na^2\hh|\| _{L^2}+C\|\bar\xx^{\delta_0/2} |\hh||\na^2\uu|\|_{L^2}\no
\leq &C\|\na\uu\|_{L^{2q/(q-2)}}\|\bar\xx^{a}\na\hh\|_{L^q} +C\|\bar\xx^{-\delta_0/2} \uu\|_{L^\infty} \|\bar\xx^{\delta_0} \na^2\hh\|_{L^2}\no
&+C\|\bar\xx^{\delta_0/2}\hh\|_{L^\infty} \|\na^2\uu\|_{L^2}\no
\leq &C,
\end{align*}
due to \eqref{4.1}, \eqref{4.3}, \eqref{4.4} and \eqref{4.7}. Then \eqref{4.49} and integration by parts lead to
\begin{align}\la{4.50}
|I_7|+|I_8|+|I_9|+|I_{10}|\leq &C\|\hh_t\|_{L^4}^2\|\na\uu_{tt}\|_{L^2} +C\|\hh\|_{L^\infty}\|\hh_{tt}\|_{L^2}\|\na\uu_{tt}\|_{L^2}\\
\leq &C\left(\|\na\uu_t\|_{L^2}+\|\hh_t\|_{L^2}^2 +\|\na\hh_t\|_{L^2}^2\right) \|\na\uu_{tt}\|_{L^2}\no
\leq &\ve\|\na\uu_{tt}\|_{L^2}^2 +C(\ve)\left(1+\|\na\uu_t\|_{L^2}^2\right),\nonumber
\end{align}
in terms of \eqref{4.5}, \eqref{4.7}, \eqref{4.27} and \eqref{4.49}.

Substituting \eqref{4.42}, \eqref{4.43}, \eqref{4.45}--\eqref{4.48} and \eqref{4.50} into \eqref{4.41}, choosing $\ve$ suitably small, and multiplying the resultant inequality by $t^2$, we get \eqref{4.37} after using Gronwall inequality and \eqref{4.15}. Therefore, the proof of Lemma \ref{l4.4} is completed. \hfill $\Box$

\section{Proof of Theorems \ref{thm1}--\ref{thm3}}

With all the a priori estimates obtained in Section 3 and 4 at hand, now we are ready to prove the main results of this paper in this section.

{\it Proof of Theorem \ref{thm1}.} Let $(\n_0, \uu_0, \hh_0)$ be as in Theorem \ref{thm1}. Without loss of generality, we assume that the initial density $\n_0$ satisfies
$$
\int_{\mr^2}\n_0 d\xx=1,
$$
which means that there has a positive constant $N_0$ such that
\be\la{5.1}
\int_{B_{N_0}}\n_0d\xx\geq\frac34\int_{\mr^2}\n_0 d\xx=\frac34.
\ee
We construct that $\n_0^R=\hat\n_0^R+R^{-1}e^{-|\xx|^2}$ where $0\leq \hat\n_0^R\in C_0^\infty(\mr^2)$ satisfies
\be\la{5.2}
\left\{
\begin{array}{lll}
\int_{B_{N_0}}\hat\n_0^R d\xx\geq 1/2,\\
\bar\xx^a\hat\n_0^R\rightarrow\bar\xx^a\n_0 \quad {\rm in\ } L^1(\mr^2)\cap H^1(\mr^2)\cap W^{1, q}(\mr^2),
\end{array}
\right.\quad{\rm as}\ R\rightarrow\infty.
\ee
Then, we choose $\hh_0^R\in \{\mathbf{w}\in C_0^\infty(B_R)|\div \mathbf{w}=0\}$ satisfying
\be\la{5.3}
\hh_0^R\bar\xx^a\rightarrow \hh_0\bar\xx^a \quad  {\rm in \ }H^1(\mr^2)\cap W^{1, q}(\mr^2),\qquad {\rm as}\ R \rightarrow\infty.
\ee
Next, since $\na\uu_0\in L^2(\mr^2)$, choosing $\mathbf v_i^R\in C_0^\infty(B_R)$ $(i=1, 2)$ such that
\be\la{5.4}
\lim_{R\rightarrow\infty}\|\mathbf v_i^R-\p_i\uu_0\|_{L^2(\mr^2)}=0,\quad i=1, 2,
\ee
we consider the unique smooth solution $\uu_0^R$ of the following elliptic problem:
\be\la{5.5}
\left\{
\begin{array}{lll}
-\triangle \uu_0^R+\nn\uu_0^R=\sqrt{\nn}\mathbf{h}^R-\p_i\mathbf v_i^R, &\ {\rm in\ }B_R,\\
\uu_0^R=0, &\ {\rm on\ }\p B_R,
\end{array}
\right.
\ee
where $\mathbf{h}^R=(\sqrt{\n_0}\uu_0)*j_{1/R}$ with $j_\delta$ being the standard mollifying kernel of width $\delta$. Extending $\uu_0^R$ to $\mr^2$ by defining $0$ outside of $B_R$ and denoting it by $\ti\uu_0^R$, we claim that
\be\la{5.6}
\lim_{R\rightarrow\infty}\left(\|\na(\ti\uu_0^R-\uu_0)\|_{L^2(\mr^2)} +\|\sqrt{\nn}\ti\uu_0^R-\sqrt{\n_0}\uu_0\|_{L^2(\mr^2)} \right)=0.
\ee

In fact, it is easy to find that $\ti\uu_0^R$ is also a solution of \eqref{5.5} in $\mr^2$. Multiplying \eqref{5.5} by $\ti\uu_0^R$ and integrating the resultant equality over $B_R$ lead to
\begin{align*}
&\int_{B_R}\nn|\ti\uu_0^R|^2d\xx+\int_{\O}|\na\ti\uu_0^R|^2d\xx\no
\leq &C\|\sqrt{\nn}\ti\uu_0^R\|_{L^2(\O)}\|\mathbf{h}^R\|_{L^2(\O)} +C\|\mathbf{v}_i^R\|_{L^2(\O)}\|\p_i\ti\uu_0^R\|_{L^2(\O)}\no
\leq & \ve\|\na\ti\uu_0^R\|_{L^2(\O)}^2+\ve\|\sqrt{\nn}\ti\uu_0^R\| _{L^2(\O)}^2 +C(\ve),
\end{align*}
which yields
\be\la{5.7}
\|\sqrt{\nn}\ti\uu_0^R\| _{L^2(\O)}^2+\|\na\ti\uu_0^R\|_{L^2(\O)}^2\leq C,
\ee
for some constant $C$ independent of $R$.

Therefore, we conclude from \eqref{5.2} and \eqref{5.7} that there exists a subsequence $R_j\rightarrow\infty$ and a function $\ti\uu_0 \in H^1_{\rm loc}(\mr^2)|\sqrt{\n_0}\ti\uu_0\in L^2(\mr^2), \na \ti\uu_0\in L^2(\mr^2)$ such that
\be\la{5.8}
\left\{
\begin{array}{lll}
\sqrt{\n_0^{R_j}}\ti\uu_0^{R_j}\rightharpoonup \sqrt{\n_0}\ti\uu_0 &\ {\rm weakly\ in\ }L^2(\mr^2),\\
\na \ti\uu_0^{R_j}\rightharpoonup \na\ti\uu_0 &\ {\rm weakly\ in\ }L^2(\mr^2).
\end{array}
\right.
\ee
Next, it follows from \eqref{5.4}, \eqref{5.5} and \eqref{5.8} that, for any $\psi\in C_0^\infty(\mr^2)$,
$$
\int_{\mr^2}\p_i(\ti\uu_0-\uu_0)\cdot\p_i\psi d\xx+\int_{\mr^2} \n_0(\ti\uu_0-\uu_0)\cdot\psi d\xx=0,
$$
which gives that
\be\la{5.10}
\ti\uu_0=\uu_0.
\ee
Furthermore, we get from \eqref{5.5} that
$$
\underset{R_j\rightarrow\infty}{\lim\sup} \int_{\mr^2}\left(|\na \ti\uu_0^{R_j}|^2+\n_0^{R_j}|\ti\uu_0^{R_j}|^2\right)d\xx \leq \int_{\mr^2}\left(|\na\uu_0|^2+\n_0|\uu_0|^2\right)d\xx,
$$
which combined with \eqref{5.8} shows
$$
\lim_{R_j\rightarrow\infty}\int_{\mr^2}|\na \ti\uu_0^{R_j}|^2d\xx= \int_{\mr^2}|\na\uu_0|^2d\xx,\quad\lim_{R_j\rightarrow\infty}\int_{\mr^2} \n_0^{R_j}|\ti\uu_0^{R_j}|^2 d\xx=\int_{\mr^2}\n_0|\uu_0|^2d\xx.
$$
This, along with \eqref{5.8} and \eqref{5.10}, yields \eqref{5.6}.

Then, due to Lemma \ref{l2.1}, the initial-boundary value problem \eqref{2.2} with the initial data $(\n_0^R, \uu_0^R, \hh_0^R)$ has a classical solution $(\n^R, \uu^R, \hh^R)$ on $B_R\times[0, T_R]$. Moreover, Proposition \ref{prop} gives that there has been a $T_0$ independent of $R$ such that \eqref{3.3} holds for $(\n^R, \uu^R, \hh^R)$. Extending $(\n^R, \uu^R, \hh^R)$ by zero on $\mr^2\backslash B_R$ and denoting it by
$$
\ti\n^R\triangleq\varphi_R\n^R, \quad \ti\uu^R, \quad \ti\hh^R\triangleq\varphi_R\hh^R,
$$
with $\varphi_R$ as in \eqref{3.6}, we first deduce from \eqref{3.3} that
\be\la{5.11}
\sup_{0\leq t\leq T_0}\left(\|\sqrt{\ti\n^R}\ti\uu^R\|_{L^2} +\|\na \ti\uu^R\|_{L^2}\right)
\leq C+C\sup_{0\leq t\leq T_0}\|\na\uu^R\|_{L^2}
\leq C,
\ee
and
\be\la{5.12}
\sup_{0\leq t\leq T_0}\|\ti\n^R\bar\xx^a\|_{L^1\cap L^\infty}\leq C.
\ee

Next, for $p\in [2, q]$, it follows from \eqref{3.3} and \eqref{3.33} that
\begin{align}\la{5.13}
&\sup_{0\leq t\leq T_0}\left(\|\na(\bar\xx^a\ti\n^R)\|_{L^p(\mr^2)}+ \|\bar\xx^a\na\ti\n^R\|_{L^p(\mr^2)}+\|\na(\bar\xx^a\ti\hh^R)\| _{L^p(\mr^2)} +\|\bar\xx^a\na\ti\hh^R\|_{L^p(\mr^2)}\right)\no
\leq &C\sup_{0\leq t\leq T_0}\left(\|\bar\xx^a\na\n^R\|_{L^p(\O)} +\|\bar\xx^a\n^R\na\varphi_R\|_{L^p(\O)}+\|\n^R\na\bar\xx^a\|_{L^p(\O)} \right)\\
&+C\sup_{0\leq t\leq T_0}\left(\|\bar\xx^a\na\hh^R\|_{L^p(\O)} +\|\bar\xx^a\hh^R\na\varphi_R\|_{L^p(\O)}+\|\hh^R\na\bar\xx^a\|_{L^p(\O)} \right)\no
\leq &C +C\|\bar\xx^a\n^R\|_{L^p(\O)}+C\|\bar\xx^a\hh^R\|_{L^p(\O)}\no
\leq &C.\nonumber
\end{align}

Then, it follows from \eqref{3.3} and \eqref{3.35} that
\be\la{5.14}
\int_0^{T_0}\left(\|\na^2 \ti\uu^R\|_{L^q(\mr^2)}^{(q+1)/q}+t \|\na^2\ti\uu^R\|_{L^q(\mr^2)}^2+\|\na^2\ti\uu^R\|_{L^2(\mr^2)} ^2\right)dt\leq C,
\ee
and that for $p\in [2, q]$,
\begin{align}\la{5.15}
&\int_0^{T_0}\left(\|\bar\xx\ti\n_t^R\|_{L^p(\mr^2)}^2 +\|\bar\xx\ti\hh_t^R\|_{L^p(\mr^2)}^2\right)dt\\
\leq & C\int_0^{T_0}\left(\||\bar\xx||\uu^R||\na\n^R|\|_{L^p(B_R)}^2 +\|\bar\xx\n^R\div\uu^R\|_{L^p(B_R)}^2\right)dt\no
&+C\int_0^{T_0}\left(\||\bar\xx||\uu^R||\na\hh^R|\|_{L^p(B_R)}^2 +\|\bar\xx|\hh^R||\na\uu^R|\|_{L^p(B_R)}^2\right)dt\no
\leq &C\int_0^{T_0}\|\bar\xx^{1-a}\uu^R\|_{L^\infty}^2\left( \|\bar\xx^a\na\n^R\|_{L^p(B_R)}^2 +\|\bar\xx^a\na\hh^R\|_{L^p(B_R)}^2\right)dt\no
\leq &C.\nonumber
\end{align}
Next, one derives from \eqref{3.3} and \eqref{3.24} that
\begin{align}\la{5.16}
&\sup_{0\leq t\leq T_0}t\|\sqrt{\ti\n^R}\ti\uu_t^R\|_{L^2(\mr^2)}^2 +\int_0^{T_0}t\|\na\ti\uu_t^R\|_{L^2(\mr^2)}^2dt\\
\leq& C+C\int_0^{T_0}t\|\na\uu_t^R\|_{L^2(\O)}^2 dt
\leq C.\nonumber
\end{align}

With all these estimates \eqref{5.11}--\eqref{5.16} at hand, we find that the sequence $(\ti\n^R, \uu^R, \ti{\mathbf{\hh}}^R)$ converges, up to the extraction of subsequences, to some limit $(\n, \uu, \hh)$ in the obvious weak sense, that is, as $R\rightarrow\infty$, we have
\be\la{5.17}
\bar\xx\ti{\n^R}\rightarrow \bar\xx\n, \quad \bar\xx\ti{\hh}^R\rightarrow \bar\xx\hh, \quad {\rm in\ }C(\overline{B_N}\times[0, T_0]), \ {\rm for\ any\ } N>0,
\ee
\be\la{5.18}
\bar\xx^a\ti{\n}^R\rightharpoonup \bar\xx^a\n,\quad \bar\xx^a\ti{\hh}^R\rightharpoonup \bar\xx^a\hh, \quad{\rm weakly\ *\ in\ } L^\infty(0, T_0; H^1(\mr^2)\cap W^{1, q}(\mr^2)),
\ee
\be\la{5.19}
\sqrt{\ti\n^R}\ti\uu^R\rightharpoonup\sqrt\n\uu,\quad \na\ti\uu^R\rightharpoonup\na\uu,\quad {\rm weakly\ *\ in\ } L^\infty(0, T_0; L^2(\mr^2)),
\ee
\be\la{5.20}
\na^2\ti\uu^R\rightharpoonup\na^2\uu,\quad{\rm weakly\ in\ } L^{(q+1)/q}(0, T_0; L^q(\mr^2))\cap L^2((0, T_0)\times\mr^2),
\ee
\be\la{5.21}
t^{1/2}\na^2\ti\uu^R\rightharpoonup t^{1/2}\na^2\uu,\quad{\rm weakly\ in\ } L^{2}(0, T_0; L^q(\mr^2)),
\ee
\be\la{5.22}
t^{1/2}\sqrt{\ti\n^R}\ti\uu^R_t\rightharpoonup t^{1/2}\sqrt\n\uu_t, \quad \na\ti\uu^R\rightharpoonup\na\uu,\quad{\rm weakly\ *\ in\ } L^{\infty}(0, T_0; L^2(\mr^2)),
\ee
\be\la{5.23}
t^{1/2}\na\ti\uu^R_t\rightharpoonup t^{1/2}\na\uu_t,\quad{\rm weakly\ in\ } L^{\infty}((0, T_0)\times\mr^2),
\ee
and
\be\la{5.24}
\bar\xx^a\n\in L^\infty(0, T_0; L^1(\mr^2)),\quad\inf_{0\leq t\leq T_0} \int_{B_{2N_0}}\n(t, \xx)d\xx\geq \frac14.
\ee

Then, letting $R\rightarrow\infty$, it follows from \eqref{5.17}--\eqref{5.24} that $(\n, \uu, \hh)$ is a strong solution of \eqref{1.1}--\eqref{1.7} on $(0, T_0]\times\mr^2$ satisfying \eqref{1.10} and \eqref{1.11}. Therefore, the proof of the existence part of Theorem \ref{thm1} is completed.

It solely remains to prove the uniqueness of the strong solution satisfying \eqref{1.10} and \eqref{1.11}. Let $(\n_1, \uu_1, \hh_1)$
and $(\n_2, \uu_2, \hh_2)$ be two strong solutions satisfying \eqref{1.10} and \eqref{1.11} with the same initial data, and denote
$$
\Psi\triangleq\n_1-\n_2,\quad U\triangleq\uu_1-\uu_2,\quad \Phi\triangleq \hh_1-\hh_2.
$$

First, subtracting the mass equation satisfied by $(\n_1, \uu_1, \hh_1)$
and $(\n_2, \uu_2, \hh_2)$ yields that
\be\la{5.25}
\Psi_t+\uu_2\cdot\na\Psi+\Psi\div\uu_2+\n_1\div U+U\cdot\na\n_1=0.
\ee
Multiplying \eqref{5.25} by $2\Psi\bar\xx^{2r}$ for $r\in (1, \ti a)$ with $\ti a=\min\{2, a\}$, and integrating by parts give
\begin{align*}
\frac{d}{dt}\|\Psi\bar\xx^r\|_{L^2}^2\leq &C\left(\|\na\uu_2\|_{L^\infty} +\|\uu_2\bar\xx^{-1/2}\|_{L^\infty}\right)\|\Psi\bar\xx^r\|_{L^2}^2 +C\|\n_1\bar\xx^r\|_{L^\infty}\|\na U\|_{L^2}\|\Psi\bar\xx^r\|_{L^2}\no
&+C\|\Psi\bar\xx^r\|_{L^2}\|U\bar\xx^{-(a-r)}\|_{L^{2q/[(q-2)(\ti a-r)]}} \|\bar\xx^a\na\n_1\|_{L^{2q/[q-(q-2)(\ti a-r)]}}\no
\leq &C\left(1+\|\na\uu_2\|_{W^{1, q}}\right) \|\Psi\bar\xx^r\|_{L^2}^2 +C\|\Psi\bar\xx^r\|_{L^2}\left(\|\sqrt{\n_1} U\|_{L^2}+\|\na U\|_{L^2}\right),
\end{align*}
which together with Gronwall inequality yields that for all $0\leq t\leq T_0$
\be\la{5.26}
\|\Psi\bar\xx^r\|_{L^2}\leq C\int_0^t\left(\|\sqrt{\n_1} U\|_{L^2}+\|\na U\|_{L^2}\right)ds.
\ee

Then, subtracting the magnetic equation satisfied by $(\n_1, \uu_1, \hh_1)$ and $(\n_2, \uu_2, \hh_2)$ leads to
\be\la{5.27}
\Phi_t+\uu_1\cdot\na\Phi+U\cdot\na\hh_2+\hh_1\div U+\Phi \div \uu_2=\hh_1\cdot\na U+\Phi\cdot\na \uu_2.
\ee
Multiplying \eqref{5.27} by $2\Phi\bar\xx^{2r}$ and integrating by parts show that
\begin{align*}
&\frac{d}{dt}\|\Phi\bar\xx^r\|_{L^2}^2\no
\leq &C\left(\|\na\uu_2\|_{L^\infty} +\|\na\uu_1\|_{L^\infty}+\|\uu_1\bar\xx^{-1/2}\|_{L^\infty}\right) \|\Phi\bar\xx^r\|_{L^2}^2 +C\|\hh_1\bar\xx^r\|_{L^\infty}\|\na U\|_{L^2}\|\Psi\bar\xx^r\|_{L^2}\no
&+C\|\Phi\bar\xx^r\|_{L^2}\|U\bar\xx^{-(a-r)}\|_{L^{2q/[(q-2)(\ti a-r)]}} \|\bar\xx^{\ti a}\na\hh_2\|_{L^{2q/[q-(q-2)(\ti a-r)]}}\no
\leq &C\left(1+\|\na\uu_1\|_{W^{1, q}}+\|\na\uu_2\|_{W^{1, q}}\right) \|\Phi\bar\xx^r\|_{L^2}^2 +C\|\Phi\bar\xx^r\|_{L^2}\left(\|\sqrt{\n_1} U\|_{L^2}+\|\na U\|_{L^2}\right),
\end{align*}
which together with Gronwall inequality gives that
\be\la{5.28}
\|\Phi\bar\xx^r\|_{L^2}\leq C\int_0^t\left(\|\sqrt{\n_1} U\|_{L^2}+\|\na U\|_{L^2}\right)ds,
\ee
for all $t\in [0, T_0]$.

Next, subtracting the momentum equation satisfied by $(\n_1, \uu_1, \hh_1)$ and $(\n_2, \uu_2, \hh_2)$ shows that
\begin{align}\la{5.29}
&\n_1U_t+\n_1\uu_1\cdot\na U-\mu\triangle U-(\mu+\lambda)\na \div U\\
=&-\n_1 U\cdot\na\uu_2-\Psi(\uu_{2t}+\uu_2\cdot\na\uu_2) -\na \left(P(\n_1)-P(\n_2)\right)\no
&-\frac12\na \left(|\hh_1|^2-|\hh_2|^2\right)+\hh_1\cdot\na\Phi +\Phi\cdot\na\hh_2.\nonumber
\end{align}
Multiplying \eqref{5.29} by $U$ and integrating by parts yields that
\begin{align}\la{5.30}
&\frac{d}{dt}\|\sqrt{\n_1}U\|_{L^2}^2+2\mu\|\na U\|_{L^2}^2 +2(\mu+\lambda) \|\div U\|_{L^2}^2\\
\leq &C\|\na\uu_2\|_{L^\infty}\|\sqrt{\n_1}U\|_{L^2}^2+C\left( \|\hh_1\bar\xx^a\|_{L^\infty}+\|\hh_2\bar\xx^a\| _{L^\infty}\right)\|\Phi\bar\xx^r\|_{L^2}\|\na U\|_{L^2}\no
&+C\int |\Psi||U|(|\uu_{2t}|+|\uu_2||\na\uu_2|)d\xx\no
&+C\left(\|P(\n_1)-P(\n_2)\|_{L^2}+\||\hh_1|^2-|\hh_2|^2\|_{L^2}\right) \|\div U\|_{L^2}\no
\leq &\ve\|\na U\|_{L^2}^2+C\left(1+\|\na\uu_2\|_{L^\infty} \right) \left(\|\sqrt{\n_1}U\|_{L^2}^2+\|\Phi\bar\xx^r\|_{L^2}^2\right) +J_1+J_2.\nonumber
\end{align}
Then, H\"{o}lder inequality yields that
\begin{align}\la{5.31}
J_1\leq &C\|\Psi\bar\xx^r\|_{L^2}\|\bar\xx^{-r/2} U\|_{L^4}\left( \|\uu_{2t}\bar\xx^{-r/2}\|_{L^4}+\|\na\uu_2\|_{L^\infty}\|\uu_2 \bar\xx^{-r/2}\|_{L^4}\right)\\
\leq &C(\ve)\left(\|\sqrt{\n_2}\uu_{2t}\|_{L^2}^2+\|\na\uu_{2t}\| _{L^2}^2+\|\na\uu_2\|_{L^\infty}^2 \right)\|\Psi\bar\xx^r\|_{L^2}^2 +\ve\left(\|\sqrt{\n_1}U\|_{L^2}^2+\|\na U\|_{L^2}^2\right).\nonumber
\end{align}
Next, Lagrange's mean value theorem together with \eqref{5.26} and \eqref{5.28} gives that
\begin{align}\la{5.32}
J_2\leq &C(\|\n_1\bar\xx^a\|_{L^\infty}, \|\n_2\bar\xx^a\|_{L^\infty}+\|\hh_1\bar\xx^a\|_{L^\infty}+\|\hh_2 \bar\xx^a\|_{L^\infty})\left(\|\Psi\bar\xx^r\|_{L^2} +\|\Phi\bar\xx^r\|_{L^2}\right) \|\na U\|_{L^2}\no
\leq &\ve\|\na U\|_{L^2}^2+C(\ve)\left(\|\Psi\bar\xx^r\|_{L^2} ^2 +\|\Phi\bar\xx^r\|_{L^2}^2\right)\\
\leq &\ve\|\na U\|_{L^2}^2+C\int_0^t\left(\|\sqrt{\n_1} U\|_{L^2}^2+\|\na U\|_{L^2}^2\right)ds.\nonumber
\end{align}

Denoting
$$
G(t)\triangleq \|\sqrt{\n_1} U\|_{L^2}^2+\int_0^t\left(\|\sqrt{\n_1} U\|_{L^2}^2+\|\na U\|_{L^2}^2\right)ds,
$$
and putting \eqref{5.31} and \eqref{5.32} into \eqref{5.30} and choosing $\ve$ small enough, one gives that
$$
G'(t)\leq C\left(1+\|\na\uu_2\|_{L^\infty}+t\|\na^2\uu_2\|_{L^q}^2+ t\|\na\uu_{2t}\|_{L^2}^2\right)G(t),
$$
which together with Gronwall inequality and \eqref{1.10} yields that $G(t)=0$. Therefore, $U(t, \xx)=0$ for almost everywhere $(t, \xx)\in (0, T_0)\times\mr^2$. Then, \eqref{5.26} and \eqref{5.28} imply that $\Psi(t, \xx)=\Phi(t, \xx)=0$ for almost everywhere $(t, \xx)\in (0, T_0)\times\mr^2$. The uniqueness of the strong solution is finished and we complete the proof of Theorem \ref{thm1}. \hfill $\Box$

{\it Proof of Theorem \ref{thm2}.} Let $(\n_0, \uu_0, \hh_0)$ be as in Theorem \ref{thm2}. Without loss of generality, assume that
$$
\int_{\mr^2}\n_0d\xx=1,
$$
which implies that there exists a positive constant $N_0$ such that \eqref{5.1} holds. We construct that $\n_0^R=\hat\n_0^R+R^{-1}e^{-|\xx|^2}$ where $0\leq \hat\n_0^R\in C_0^\infty(\mr^2)$ satisfies \eqref{5.2} and
\be\la{5.33}
\left\{
\begin{array}{lll}
\na^2\hat\n_0^R\rightarrow\na^2\n_0, \quad \na^2 P(\hat\n_0^R)\rightarrow \na^2 P(\n_0), &{\rm in\ } L^q(\mr^2),\\
\bar\xx^{\delta_0}\na^2\hat\n_0^R\rightarrow \bar\xx^{\delta_0}\na^2\n_0, \quad \bar\xx^{\delta_0}\na^2P(\hat\n_0^R)\rightarrow \bar\xx^{\delta_0}\na^2P(\n_0), &{\rm in\ } L^2(\mr^2),
\end{array}
\right.
\ee
as $R\rightarrow\infty$. Then, we also choose $\hh_0^R\in \{\mathbf{w}\in C_0^\infty(B_R)|\div \mathbf{w}=0\}$ satisfying \eqref{5.3} and
\be\la{5.34}
\left\{
\begin{array}{lll}
\na^2 \hh^R_0\rightarrow\na^2\hh_0,&{\rm in\ } L^q(\mr^2),\\
\bar\xx^{\delta_0}\na^2\hh_0^R\rightarrow \bar\xx^{\delta_0}\na^2\hh_0, &{\rm in\ } L^2(\mr^2),
\end{array}
\right.\quad {\rm as\ } R\rightarrow\infty.
\ee

Then, we consider the unique smooth solution $\uu_0^R$ of the following elliptic problem
\be\la{5.35}
\left\{
\begin{array}{lll}
-\mu\triangle\uu_0^R-(\mu+\lambda)\na\div\uu_0^R+\na P(\n_0^R)\\
\qquad\qquad=(\na\times \hh_0^R)\times\hh_0^R-\n_0^R\uu_0^R+\sqrt{\n_0^R} \mathbf{h}^R, & {\rm in\ } B_R,\\
\uu_0^R=0, &{\rm on\ } \p B_R,
\end{array}
\right.
\ee
where $\mathbf{h}^R=(\sqrt{\n_0}\uu_0+\mathbf{g})*j_{1/R}$ with $j_\delta$ being the standard mollifying kernel of width $\delta$. Multiplying \eqref{5.35} by $\uu_0^R$ and integrating the resultant equation over $B_R$, it is easy to show that
\begin{align*}
&\|\sqrt{\n_0^R}\uu_0^R\|_{L^2(B_R)}^2+\mu\|\na\uu_0^R\|_{L^2(B_R)}^2 +(\mu+\lambda)\|\div\uu_0^R\|_{L^2(B_R)}^2 \no
\leq &\int_{\O}P(\n_0^R)|\div\uu_0^R| d\xx+\frac12\int_{\O} |\hh_0^R|^2|\div\uu_0^R| d\xx+\int_{\O}|\hh_0^R||\na \hh_0^R||\uu_0^R|d\xx\no
&+\|\sqrt{\n_0^R}\uu_0^R\|_{L^2(B_R)}\|\mathbf{h}^R\|_{L^2(B_R)}\no
\leq &\ve\left(\|\sqrt{\n_0^R}\uu_0^R\|_{L^2(B_R)}^2+\|\na\uu_0^R \|_{L^2(B_R)}^2 \right)+C(\ve),
\end{align*}
which implies that
\be\la{5.36}
\|\sqrt{\n_0^R}\uu_0^R\|_{L^2(B_R)}^2+\|\na\uu_0^R \|_{L^2(B_R)}^2 \leq C,
\ee
for some constant $C$ independent of $R$. Due to \eqref{2.9}, we have
\begin{align}\la{5.37}
\|\na^2\uu_0^R\|_{L^2(B_R)}\leq & C\|\na P(\n_0^R)\|_{L^2(\O)}+C\||\hh_0^R||\na\hh_0^R|\|_{L^2(\O)}\\
&+C\|\n_0^R\uu_0^R\|_{L^2(\O)}+C\|\sqrt{\n_0^R}\mathbf{h}^R \|_{L^2(\O)}\no
\leq &C.\nonumber
\end{align}

Next, extending $\uu_0^R$ to $\mr^2$ by defining $0$ outside $\O$ and denoting it by $\ti\uu_0^R$, we deduce from \eqref{5.36} and \eqref{5.37} that
$$
\|\na\ti\uu_0^R\|_{H^1(\mr^2)}\leq C,
$$
which together with \eqref{5.33} and \eqref{5.36} gives that there exists a subsequence $R_j\rightarrow\infty$ and a function $\ti\uu_0\in \{\ti\uu_0\in H^2_{\rm loc}(\mr^2)|\sqrt{\n_0}\ti\uu_0 \in L^2(\mr^2), \na\ti\uu_0\in H^1(\mr^2)\}$ such that
\be\la{5.38}
\left\{
\begin{array}{lll}
\sqrt{\n_0^{R_j}}\ti\uu_0^{R_j}\rightharpoonup \sqrt{\n_0}\ti\uu_0,&{\rm weakly\ in\ } L^2(\mr^2),\\
\na\ti\uu_0^{R_j}\rightharpoonup \na\ti\uu_0, &{\rm weakly\ in\ } H^1(\mr^2).
\end{array}
\right.
\ee
It is easy to check that $\ti\uu_0^R$ satisfies \eqref{5.35}, then one can deduce from \eqref{5.33}, \eqref{5.34}, \eqref{5.35} and \eqref{5.38} that $\ti\uu_0$ satisfies
$$
-\mu\triangle\ti\uu_0-(\mu+\lambda)\na\div\ti\uu_0+\na P(\n_0)+\n_0\ti\uu_0=(\na\times\hh_0)\times\hh_0+\n_0\uu_0 +\sqrt{\n_0}\mathbf{g},
$$
which combined with \eqref{1.13} yields that
\be\la{5.39}
\ti\uu_0=\uu_0.
\ee

Next, we get from \eqref{5.35} that
$$
\underset{R_j\rightarrow\infty}{\lim\sup} \int_{\mr^2}\left(|\na \ti\uu_0^{R_j}|^2+\n_0^{R_j}|\ti\uu_0^{R_j}|^2\right)d\xx \leq \int_{\mr^2}\left(|\na\uu_0|^2+\n_0|\uu_0|^2\right)d\xx,
$$
which combined with \eqref{5.38} shows
$$
\lim_{R_j\rightarrow\infty}\int_{\mr^2}|\na \ti\uu_0^{R_j}|^2d\xx= \int_{\mr^2}|\na\uu_0|^2d\xx,\quad\lim_{R_j\rightarrow\infty}\int_{\mr^2} \n_0^{R_j}|\ti\uu_0^{R_j}|^2 d\xx=\int_{\mr^2}\n_0|\uu_0|^2d\xx.
$$
This, along with \eqref{5.38} and \eqref{5.39}, shows that
\be\la{5.40}
\lim_{R\rightarrow\infty}\left(\|\na(\ti\uu_0^R-\uu_0)\|_{L^2(\mr^2)}+ \|\sqrt{\n_0^R}\ti\uu_0^R-\sqrt{\n_0}\uu_0\|_{L^2(\mr^2)}\right)=0.
\ee
Similar to \eqref{5.40}, we can also obtain that
$$
\lim_{R\rightarrow\infty}\|\na^2(\ti\uu_0^R-\uu_0)\|_{L^2(\mr^2)}=0.
$$

Finally, in terms of Lemma \ref{l2.1}, the initial-boundary value problem \eqref{2.2} with the initial data $(\n_0^R, \uu_0^R, \hh_0^R)$ has a classical solution $(\n^R, \uu^R, \hh^R)$ on $[0, T_R]\times\O$. Hence, there has a generic positive constant $C$ independent of $R$ such that all those estimates stated in Proposition \ref{prop} and Lemma \ref{l4.1}--\ref{l4.4} hold for $(\n^R, \uu^R, \hh^R)$. Extending $(\n^R, \uu^R, \hh^R)$ by zero on $\mr^2\setminus \O$ and denoting
$$
\ti\n^R\triangleq\varphi_R\n^R, \quad \ti\uu^R, \quad \ti\hh^R\triangleq\varphi_R\hh^R,
$$
with $\varphi_R$ as in \eqref{3.6}. We deduce from \eqref{3.3} and Lemma \ref{l4.1}--\ref{l4.4} that the sequence $(\ti\n^R, \ti\uu^R, \ti\hh^R)$ converges weakly, up to the extraction of subsequences, to some limit $(\n, \uu, \hh)$ satisfying \eqref{1.10}, \eqref{1.11} and \eqref{1.14}. Moreover, standard arguments shows that $(\n, \uu, \hh)$ is in fact a classical solution to the problem \eqref{1.1}--\eqref{1.7}. The proof of
Theorem \ref{thm2} is finished. \hfill $\Box$

{\it Proof of Theorem \ref{thm3}.} Now, we prove \eqref{1.15}. Let $(\n, \uu, \hh)$ be the unique classical solution to \eqref{1.1}--\eqref{1.7} obtained in Theorem \ref{thm2}. We assume that the opposite holds, i.e.,
\be\la{6.1}
\underset{T\rightarrow T^*}{\lim\sup}\left(\|\n\|_{L^\infty(0, T; L^\infty(\mr^2))}+\|\hh\|_{L^\infty(0, T; L^\infty(\mr^2))} +\|\na\hh\|_{L^2(0, T; L^2(\mr^2))}\right)=M_0<\infty.
\ee

First, for standard energy estimate gives
\be\la{6.2}
\sup_{0\leq t\leq T}\left(\|\sqrt\n\uu\|_{L^2}^2+\|\hh\|_{L^2}^2 +\|P(\n)\|_{L^1}\right) +\int_0^T\|\na\uu\|_{L^2}^2dt\leq C.
\ee

Based on the assumption \eqref{6.1} and the energy estimates \eqref{6.2}, we derive the following estimates \eqref{6.3}, \eqref{6.12}, \eqref{6.20} and \eqref{5.65} by Lemma \ref{l6.1}--\ref{l5.4}, respectively, which will be used to complete the proof of Theorem \ref{thm3}. In the rest of this section, the generic positive constant $C$ may depend on $M_0$, $\mu, \lambda, \ga, q, a, \eta_0, N_0$, $\|\mathbf{g}\|_{L^2}$, $C_0$, $\delta_0$ and initial data assumed in Theorem \ref{thm2}.

\bl\la{l6.1}
Let $(\n, \uu, \hh)$ be classical solution obtained in Theorem \ref{thm2}. Under the condition \eqref{6.1}, it holds that for $0\leq t\leq T^*$
\be\la{6.3}
\sup_{0\leq t\leq T}\|\na\uu\|_{L^2}^2+\int_0^{T}\|\sqrt\n\dot\uu\| _{L^2}^2 dt\leq C.
\ee
\el
\pf Multiplying \eqref{1.2} by $\dot\uu$ and integrating by parts over $\mr^2$, direct calculations yield that
\begin{align}\la{6.4}
\int\n|\dot\uu|^2d\xx =&-\int\dot\uu\cdot\na P(\n)d\xx+\mu\int\dot\uu \cdot \triangle \uu d\xx+(\mu+\lambda)\int\dot\uu\cdot\na\div\uu d\xx\\
&-\frac12\int\dot\uu\cdot\na |\hh|^2d\xx+\int\hh\cdot\na\hh \cdot\dot\uu d\xx.\nonumber
\end{align}

Now, we estimate each term on the right-hand side of \eqref{6.4}. First, it follows from \eqref{3.14} and integration by parts that
\begin{align}\la{6.5}
-\int\dot\uu\cdot\na P(\n)d\xx=&\frac{d}{dt}\int P(\n)\div\uu d\xx+ \int\left[(\ga-1)P(\n)(\div \uu)^2+P(\n)\p_i\uu_j\p_j\uu_i\right]d\xx\no
\leq &\frac{d}{dt}\int P(\n)\div\uu d\xx+ C\|\na\uu\|_{L^2}^2.
\end{align}

Then, integration by parts lead to
\begin{align}\la{6.6}
\mu\int\triangle \uu\cdot\dot\uu d\xx=&-\frac\mu2\frac{d}{dt}\int|\na\uu|^2d\xx+\mu\int\triangle \uu (\uu\cdot\na)\uu d\xx\\
\leq &-\frac\mu2\frac{d}{dt}\int|\na\uu|^2d\xx +C\|\na\uu\|_{L^3}^3, \nonumber
\end{align}
and similarly,
\be\la{6.7}
(\mu+\lambda)\int\na\div\uu\cdot\dot\uu d\xx\leq -\frac{\mu+\lambda}{2}\frac{d}{dt}\int(\div\uu)^2d\xx +C\|\na\uu\|_{L^3}^3.
\ee

Next, \eqref{1.3} and integration by parts yield that
\begin{align}\la{6.8}
-\frac12\int\dot\uu\cdot\na |\hh|^2d\xx=&\frac12\frac{d}{dt}\int |\hh|^2\div\uu d\xx-\frac12\int |\hh|^2\na\uu\div\uu d\xx\\
&-\frac12\int |\hh|^2\na\uu\cdot\na\uu d\xx+\frac12\int |\hh|^2(\div\uu)^2d\xx\no
\leq & \frac12\frac{d}{dt}\int |\hh|^2\div\uu d\xx +C\|\na\uu\|_{L^2}^2,\nonumber
\end{align}
and similarly, that
\be\la{6.9}
\int\hh\cdot\na\hh \cdot\dot\uu d\xx \leq \frac{d}{dt}\int \hh\cdot\na \uu\cdot\hh d\xx+C\|\na\uu\|_{L^2}^2.
\ee

Substituting \eqref{6.5}--\eqref{6.9} into \eqref{6.4}, it follows from \eqref{2.14} and \eqref{6.1} that
\begin{align}\la{6.10}
A'(t)+\|\sqrt\n \dot\uu\|_{L^2}^2\leq &C\|\na\uu\|_{L^2}^2 +C\|\na\uu\|_{L^3}^3\\
\leq&C\|\na\uu\|_{L^2}^2 +\frac12\|\sqrt\n \dot\uu\|_{L^2}^2 +C\|\na\uu\|_{L^2}^4+C\|\na\hh\|_{L^2}^2+C.\nonumber
\end{align}
where
\begin{align}\la{6.11}
A(t)\triangleq & \frac\mu 2\|\na\uu\|_{L^2}^2+ \frac{\mu+\lambda}{2}\|\div\uu\|_{L^2}^2 -\int\div\uu P(\n) d\xx\\ &-\frac12\int|\hh|^2\div\uu d\xx+\int \hh\cdot\na\uu\cdot\hh d\xx\no
\geq&\frac\mu4\|\na\uu\|_{L^2}^2+ \frac{\mu+\lambda}{2} \|\div\uu\|_{L^2}^2 -C, \nonumber
\end{align}
due to \eqref{6.1} and \eqref{6.2}. Then integrating \eqref{6.10} over $(0, T)$, together with \eqref{6.1}, \eqref{6.2} and \eqref{6.11}, Gronwall inequality leads to \eqref{6.3}. Therefore, we complete the proof of Lemma \ref{l6.1}. \hfill $\Box$

\bl\la{l6.2}
Let $(\n, \uu, \hh)$ be classical solution obtained in Theorem \ref{thm2}. Under the condition \eqref{6.1}, it holds that for $0\leq t\leq T^*$
\be\la{6.12}
\sup_{0\leq t\leq T}\|\sqrt\n\dot\uu\|_{L^2}^2+\int_0^{T}\|\na\dot\uu\| _{L^2}^2 dt\leq C.
\ee
\el
\pf Operating $\dot\uu^j [\frac{\partial}{\partial t} +\div(\uu\cdot)]$ on \eqref{1.2}$^j$, summing with respect to $j$, and integrating the resultant equation by parts over $\mr^2$, one obtains
\begin{align}\la{6.13}
&\frac12\frac{d}{dt}\left(\int\n|\dot\uu|^2d\xx\right)\\
=&-\int\dot\uu^j[\pa_j P_t(\n)+\div(\pa_jP(\n)\uu)]d\xx +\mu\int\dot\uu^j[\triangle \uu^j_t+\div(\uu\triangle \uu^j)]d\xx\no
&+(\mu+\lambda)\int\dot\uu^j[\pa_t\pa_j\div\uu +\div(\uu\pa_j\div\uu)]d\xx \no &-\frac12\int\dot\uu^j\left(\p_t\p_j|\hh|^2 +\div(\uu\p_j|\hh|^2)\right)d\xx\no
&+\int\dot\uu^j\left(\p_t(\hh\cdot\na\hh^j)+\div(\uu(\hh\cdot\na \hh^j))\right)d\xx.\nonumber
\end{align}

Now we estimate each term on the righthand of \eqref{6.13}. First, it follows from \eqref{3.14}, integration by parts and careful calculations, that
\begin{align}\la{6.14}
&-\int\dot\uu^j[\pa_j P_t(\n)+\div(\pa_jP(\n)\uu)]d\xx\\
=&\int\left[-P'(\n)\n\div\uu\pa_j\dot\uu^j+\pa_k(\pa_j\dot\uu^j \uu^k)P(\n)-P(\n)\pa_j(\pa_k\dot\uu^j\uu^k)\right]d\xx\no
\leq&\ve\|\na\dot\uu\|_{L^2}^2+C\|\na\uu\| _{L^2}^2.\nonumber
\end{align}

Then, due to integration by parts and after subtle calculations, we have
\begin{align}\la{6.15}
&-\mu\int\dot\uu^j[\triangle \uu^j_t+\div(\uu\triangle \uu^j)]d\xx\\
=&-\mu\int\left[|\na\dot\uu|^2+\pa_i\dot\uu^j\pa_k\uu^k\pa_i\uu^j -\pa_i\dot\uu^j\pa_i\uu^k\pa_k\uu^j-\pa_i\uu^j\pa_i\uu^k\pa_k\dot\uu^j \right]d\xx\no
\leq &-\mu\|\na\dot\uu\|_{L^2}^2+C\|\na\uu\|_{L^4}^4,\nonumber
\end{align}
and similarly, that
\be\la{6.16}
(\mu+\lambda)\int\dot\uu^j[\pa_t\pa_j\div\uu +\div(\uu\pa_j\div\uu)]d\xx \leq -(\mu+\lambda)\|\div\dot\uu\|_{L^2}^2+C\|\na\uu\|_{L^4}^4.
\ee

Next, it follows from \eqref{1.3} and integration by parts, that
\begin{align}\la{6.17}
&-\frac12\int\dot\uu^j\left(\p_t\p_j|\hh|^2 +\div(\uu\p_j|\hh|^2)\right)d\xx\\
=&\frac12\int\p_j\dot\uu^j\div\uu|\hh|^2d\xx-\frac12\int\p_j\uu\cdot\na \dot\uu^j |\hh|^2d\xx+\int\p_j\dot\uu^j\hh\cdot(\hh\cdot\na\uu- \hh\div\uu)d\xx\no
\leq &\ve\|\na\dot\uu\|_{L^2}^2+C\|\na\uu\|_{L^2}^2,\nonumber
\end{align}
and similarly calculations lead to
\begin{align}\la{6.18}
\int\dot\uu^j\left(\p_t(\hh\cdot\na\hh^j)+\div(\uu(\hh\cdot\na \hh^j))\right)d\xx\leq \ve\|\na\dot\uu\|_{L^2}^2+C\|\na\uu\|_{L^2}^2.
\end{align}

Then, inserting \eqref{6.14}--\eqref{6.18} into \eqref{6.13} and choosing $\ve$ suitably small leads to
\be\la{6.19}
\frac{d}{dt}\|\sqrt\n\dot\uu\|_{L^2}^2+\|\na\dot\uu\| _{L^2}^2\leq C\|\n\dot\uu\|_{L^2}^2+C\|\na\hh\|_{L^2}^2+C,
\ee
where we have used \eqref{2.14}, \eqref{6.1}, \eqref{6.3}. Then, it follows from \eqref{6.19} and the compatibility conditions \eqref{1.13},  we obtain \eqref{6.12} after using Gronwall inequality and \eqref{6.1}. Thus, we complete the proof of Lemma \ref{l6.2}. \hfill $\Box$

\bl\la{l6.3}
Let $(\n, \uu, \hh)$ be classical solution obtained in Theorem \ref{thm2}. Under the condition \eqref{6.1}, it holds that for $0<T\leq T^*$
\be\la{6.20}
\sup_{0\leq t\leq T}\left(\|\n\bar\xx^a\|_{L^1\cap H^1\cap W^{1, q}} +\|\hh\bar\xx^a\|_{H^1\cap W^{1, q}}\right)\leq C.
\ee
\el
\pf First, multiplying \eqref{1.1} by $\bar\xx^a$, integrating the resultant equality over $\mr^2$, lead to
\begin{align}\la{6.21}
\frac{d}{dt}\int\n\bar\xx^ad\xx\leq &C\int \n |\uu|\bar\xx^{a-1} \ln^{1+\eta_0}(e+|\xx|^2)d\xx\\
\leq &C\|\n\bar\xx^{a-1+8/(8+a)}\|_{L^{(8+a)/(7+a)}}\|\uu \bar\xx^{-4/(8+a)}\|_{L^{8+a}}\no
\leq &C,\nonumber
\end{align}
where in the last inequality we have used \eqref{2.6}, \eqref{6.2} and \eqref{6.3}.

Next, it follows from \eqref{3.36} and \eqref{a3.40} that for $p\in [2, q]$, we have
\begin{align}\la{6.22}
&\frac{d}{dt}\left(\|\ti\hh\|_{L^2}+\|\na\ti\hh\|_{L^p}+\|w\|_{L^2}+\|\na w\|_{L^p}\right)\\
\leq &C\left(1+\|\na\uu\|_{L^\infty}+\|\uu\cdot\na \ln\bar\xx\| _{L^\infty} \right)\left(\|\ti\hh\|_{L^2}+\|\na\ti\hh\|_{L^p}+\|w\|_{L^2} +\|\na w\|_{L^p}\right)\no
&+C\left(\||\na\uu||\na\ln\bar\xx|\|_{L^p}+\||\uu||\na^2 \ln\bar\xx|\|_{L^p}+\|\na^2\uu\|_{L^p}\right)\no
\leq &\left(1+\|\na\uu\|_{L^\infty}+\|\na\uu\|_{L^3}\right)\left(1+ \|\ti\hh\| _{L^2}+\|\na\ti\hh\|_{L^q}+\|w\|_{L^2}+\|\na w\|_{L^q}\right)\no
&+C\|\na\uu\|_{L^p}+C\|\na^2\uu\|_{L^p}\no
\leq & C\left(1+\|\na\uu\|_{L^\infty}+\|\na\hh\|_{L^2}^2 +\|\na\dot\uu\|_{L^2}^2\right)\left(1+ \|\ti\hh\| _{L^2}+\|\na\ti\hh\|_{L^q}+\|w\|_{L^2}+\|\na w\|_{L^q}\right).\nonumber
\end{align}
where in the second inequality we have used \eqref{6.2}, \eqref{6.3} and similar discusses as that in \eqref{3.37} and \eqref{3.38}; in the third inequality, we have used the following facts:
\begin{align*}
\|\na^2\uu\|_{L^p}\leq &C\left(\|\n\dot\uu\|_{L^p}+\|\na\hh\|_{L^p} +\|\na \n\|_{L^p}\right)\no
\leq &C\left(1+\|\na\dot\uu\|_{L^2}+\|\na\ti\hh\|_{L^p} +\|\na w\|_{L^p}\right),
\end{align*}
which come from $L^p$-estimates on the Lam\'{e} system. It follows from \eqref{2.15}, that
\begin{align}\la{5.62}
\|\na\uu\|_{L^\infty}\leq & C\left(\|\div\uu\|_{L^\infty}+\|\o\| _{L^\infty} \right)\ln\left(e+\|\na^2\uu\|_{L^q}\right)+C\|\na\uu \|_{L^2}^2+C\\
\leq &C\left(1+\|\n\dot\uu\|_{L^q}^{q/(2(q-1))} +\|\na\hh\|_{L^q} ^{q/(2(q-1))}+\|\na\n\|_{L^q} ^{q/(2(q-1))}\right)\no
&\quad\cdot\ln\left(e+\|\n\dot\uu\|_{L^q}+\|\na\ti\hh\|_{L^q} +\|\na w\|_{L^q}\right)+C\no
\leq &C\left(1+\|\na\dot\uu\|_{L^2}^2 +\|\na\ti\hh\|_{L^q}+\|\na w\|_{L^q}\right)\ln\left(e+\|\na\ti\hh\|_{L^q}+\|\na w\|_{L^q}\right) +C,\nonumber
\end{align}
which together with \eqref{6.3}, \eqref{6.12}, \eqref{6.21} and \eqref{6.22}, after using Gronwall inequality give
\be\la{5.63}
\|\ti\hh\|_{L^2}+\|\na\ti\hh\|_{L^q}+\|w\|_{L^2}+\|\na w\|_{L^q}\leq C.
\ee
Then, it follows from \eqref{6.12}, \eqref{5.62} and \eqref{5.63} that
$$
\int_0^T\|\na\uu\|_{L^\infty}dt\leq  C,
$$
which together with \eqref{6.1}, \eqref{6.12} and \eqref{6.22} leads to
\be\la{5.64}
\sup_{0\leq t\leq T}\left(\|\na\ti\hh\|_{L^2}+\|\na w\|_{L^2}\right)\leq C.
\ee
Therefore, combining \eqref{6.21}, \eqref{5.63} and \eqref{5.64}, one obtains \eqref{6.20} and the proof of Lemma \ref{l6.3} is finished. \hfill $\Box$

With the a priori estimates obtained in Lemma \ref{l6.1}--\ref{l6.3} at hand, the following higher order estimates of the solutions which are needed to guarantee the local strong solutions to be a classical one are similar to those obtained in Lemma \ref{l4.1}--\ref{l4.4}, so we omit their detailed proofs here.

\bl\la{l5.4} Let $(\n, \uu, \hh)$ be classical solution obtained in Theorem \ref{thm2}. Under the condition \eqref{6.1}, it holds that for $0\leq t\leq T^*$
\begin{align}\la{5.65}
&\sup_{0\leq t\leq T}\Big(\|\bar\xx^{\delta_0}\na^2\n\|_{L^2} +\|\bar\xx^{\delta_0}\na^2P(\n)\|_{L^2}+\|\bar\xx^{\delta_0} \na^2\hh\|_{L^2}+t\|\na\uu_t\|_{L^2}^2+\|\na^2\n\|_{L^q}\\
&\qquad+\|\na^2P(\n)\|_{L^q} +\|\na^2\hh\|_{L^q}+t\|\na^3\uu\|_{L^2\cap L^q}+t\|\na^2\uu_t\|_{H^1} +t\|\na^2(\n\uu)\|_{L^{(q+2)/2}}\Big)\no
&\qquad+\int_0^{T_0}\left( t\|\sqrt\n \uu_{tt}\|_{L^2}^2 +t\|\na^2\uu_{t} \|_{L^2}^2+t^2\|\na\uu_{tt}\| _{L^2}^2+t^2 \|\bar\xx^{-1} \uu_{tt}\|_{L^2}^2\right)dt\leq C.\nonumber
\end{align}
\el

Now, we are ready to finish the proof of Theorem \ref{thm3}. In fact, in view of the estimates obtained in Lemma \ref{l6.1}--\ref{l5.4}, one easily shows that the functions
$$(\n, \uu, \hh)(T^*, \xx)=\lim_{T\rightarrow T^*}(\n, \uu, \hh)(T, \xx)$$
satisfying the conditions imposed on the initial data \eqref{1.9} and \eqref{1.12}. Therefore, we can take $(\n, \uu, \hh)(T^*, \xx)$ as the initial data and applying Theorem \ref{thm2} to extend our local classical solution beyond $t> T^*$. This contradicts the assumption on $T^*$. Therefore, we complete the proof of Theorem \ref{thm3}. \hfill $\Box$

\section*{Acknowledge}
This work  of M.  Chen is partially supported by   National Natural Science Foundation of China (11471191), the Independent Innovation Foundation of Shandong University (No. 2013ZRQP001) and the National Science Foundation of Shandong province of China under grant (No. ZR2015AM019) . The work of A. Zang is supported in part  National Natural Science Foundation of China (11571279)   and a part of Project GJJ151036 supported by Education Department of Jiangxi Province and  partly supported by Youth Innovation Group of Applied Mathematics  in Yichun University(2012TD006).

\end{document}